\documentclass{article}
\usepackage[a4paper, margin=1in]{geometry}
\usepackage{amsmath, amssymb, amsthm}
\usepackage{graphicx}
\usepackage{hyperref}
\usepackage{amsmath}     				%obligatory package
\usepackage{amssymb}   					%obligatory package
\usepackage{amsthm}     				%obligatory package
\usepackage{graphicx}    				%obligatory package

\usepackage{textcomp}
\usepackage[utf8]{inputenc}
\usepackage[english]{babel}    				%obligatory package
\usepackage[T1]{fontenc} 				%obligatory package
\usepackage{marvosym}    				%obligatory package (for the \Letter symbol for indicating corresponding author)
\usepackage[sc]{mathpazo}  	 			%obligatory package (font package)

\usepackage{authblk} 					%obligatory package
\usepackage[usenames]{xcolor}  			%obligatory package
\usepackage{lastpage} 					%obligatory package

\usepackage{enumerate}	 				%recommended package

\usepackage{hyperref} 					%Use this in pdflatex mode  %obligatory package

\usepackage{amssymb}
\newcommand{\N}{{\mathbb N}}
\newcommand{\cF}{{\mathcal F}}
\newcommand{\cL}{{\mathcal L}}
\newcommand{\cT}{{\mathcal T}}

\newcommand{\fF}{{\mathfrak F}}
\newcommand{\norm}[1]{\left\|#1\right\|}				% norm
\renewcommand{\d}{\,{\mathrm d}}
\newcommand{\tm}{\times}
\newcommand{\fall}{\quad\text{for all }}
\newtheorem{theorem}{Theorem}[section]

\newtheorem{lemma}[theorem]{Lemma}

\theoremstyle{definition}
\newtheorem{definition}[theorem]{Definition}
\theoremstyle{definition}
\newtheorem{remark}[theorem]{Remark}
\theoremstyle{definition}

\numberwithin{equation}{section}
\numberwithin{table}{section}
\numberwithin{figure}{section}
\title{Smooth linearization of contractive random dynamical systems in continuous time}
\author{Iryna Vasylieva \\
        \small University of Klagenfurt, A--9020 Universit\"atsstraße 65--67, Klagenfurt,  Austria \\
        \small \texttt{iryna.vasylieva@aau.at}}
\date{}

\begin{document}

\maketitle

\begin{abstract}
    We establish that uniformly exponentially stable random dynamical systems on the half line have equivalent dynamics through a  $C^m-$ conjugacy.
This result was obtained for  random  differential equations as well as for  random dynamical systems with a uniformly exponentially stable linear part.
\end{abstract}
\smallskip
\section{Introduction} % Please enter the title of your first section (only the first letter of the title should be capital).

The Hartman-Grobman theorem is a fundamental  result in the local theory of dynamical systems. It establishes a connection between the qualitative behavior of a nonlinear system and its linearization around an equilibrium point. The classical theorem states that if a nonlinear system possesses a hyperbolic equilibrium point, then the qualitative behavior of the system near the equilibrium point is topologically equivalent to that of its linearization.

The original theorem was obtained by Philip Hartman \cite{hartman:60a}, \cite{hartman:60b} and David Grobman \cite{grobman:59} independently in the 1960s. It has since become a cornerstone of dynamical systems theory. This theorem provides a powerful tool for analyzing, e.g., stability of nonlinear systems by reducing the problem to the study of linear systems.

It should be mentioned that the standard framework in Hartman-Grobman theory relies on a fundamental assumption about the hyperbolicity of the linearized system. However, when eigenvalues lie on the imaginary axis, the dynamics can undergo significant changes. For this reason, the original paper on the linearization theorem has inspired pioneering work and extensive research, leading to numerous generalizations and extensions.

Speaking about generalizations, it is important to mention results for autonomous systems in Banach spaces \cite{Pal_68}, \cite{Pugh_69},  linearization for nonautonomous problems \cite{lu2024}, \cite{palmer:73}, \cite{poetzsche:07b}  and results for Carath{\'e}odory-type differential equations in Banach spaces \cite{aulbach:wanner:00}.

In recent years, there has been growing interest in extending the Hartman-Grobman theorem to the case of random dynamical systems. These systems incorporate the effects of stochasticity or randomness into the evolution of dynamical systems (see, e.g., \cite{Arnold}, \cite{Chu_02}), and are particularly relevant in modeling real-world phenomena where uncertainties or random perturbations play a significant role.
The generalization of the Hartman-Grobman theorem to random dynamical systems aims to establish similar qualitative connections between the behavior of a nonlinear system and its linearization in the presence of random perturbations. The goal is to understand how the stability properties, attractors, and bifurcations of the deterministic system are affected by the stochastic components. Accordingly, extensions of the Hartman–Grobman theorem to random dynamical systems go back to the works  by Arnold~\cite{Arnold}  resulting from Wanner \cite{Wanner_95}.
As in the deterministic case, there are many generalizations of the random linearization theorem, for instance, in  \cite{Zhao_20}, \cite{Li_05}, \cite{Lu_05}, and others.

In linearization theory, a compelling research question revolves around investigating the smoothness properties of the linearization map, since this allows a better understanding of the information that is carried from one system to another. Smooth linearization ensures that invariant manifolds and their geometric structures, such as stable and unstable tangent directions, are consistently preserved. This is important for accurately capturing the flow near fixed points, periodic orbits, or other invariant sets, allowing  to study complex dynamics in a simplified framework.
This research problem is widely explored in deterministic settings.
For example, in \cite{shi:95}, \cite{Ray_98}, \cite{Barreira_07}  H\"older properties of linearization homeomorphisms were established.
The $C^1-$ linearization result
for contractions has been obtained in \cite{Rodrigues_17}, \cite{Zhang_11}, \cite{Back_22}.
Significant contributions in the  higher order smoothness  within the autonomous context are the results of   Sternberg \cite{Sternberg}, \cite{Sternberg_1}, who proved that under certain non-resonance conditions the
local conjugation between two systems can actually be a $C^m$  map. Later,
Belickii \cite{Bel}, \cite{Bel_1} achieved the same result under assumptions related to the spectrum
of the diffeomorphism that defines the dynamics.

The authors in \cite{Castaneda_15} delved into the differentiability properties within a particular case of Palmer's linearization theorem \cite{palmer:73}, particularly when restricting the problem to the nonnegative real line. Within this framework, sufficient conditions ensuring that the conjugating map is a \( C^2 \)-preserving orientation diffeomorphism were obtained. These conditions were established under the premise that the linear system was uniformly asymptotically stable. Recent generalizations and improvements of this result were presented in \cite{crm:20}.

Like in deterministic scenarios, the smoothness properties of the conjugating map in random dynamical systems were equally intriguing. At this point, it is worth mentioning the paper \cite{backes2021}, where conditions for Hölder continuity of the conjugacies on bounded subsets for infinite-dimensional random dynamical systems with a non-invertible linear part were established.

In the paper \cite{Valls_06}, the Hölder continuity of conjugacy for both maps and flows was established. Notably, the analysis demonstrated that such continuity could be achieved without imposing any nonresonance condition. In \cite{Lu_20}, a proof of a \( C^{1,\beta} \) generalization of the random Hartman-Grobman theorem was presented. Additionally, temperedness for the conjugacy was obtained, further enhancing the understanding of the relationship between random dynamical systems and their associated linearizations.

In \cite{Duan_14}, stochastic differential equations (SDEs) with non-Gaussian Lévy processes for two-sided time were studied. Topological equivalence between these SDEs and transformed random differential equations (RDEs) with Lévy processes was established. This allowed for a stochastic Hartman-Grobman theorem for the linearization of SDEs to be proved.

The aim of the present paper is to derive sufficient conditions for
\( C^m \)-linearizations of random dynamical systems on the nonnegative half-line. The advantage of considering the positive half-line is that no information about the past is required, allowing the focus to be solely on future dynamics. Specifically, the main contribution of the paper can be considered a generalization of the approach proposed in \cite{crm:20} in several directions.

Firstly, in the class of Carathéodory differential equations, we not only derived sufficient conditions for the topological equivalence of the nonlinear equation and its linearization, but also provided additional information about the conjugating map. Specifically, we established that the conjugating map is globally Lipschitz and belongs to the class \( C^m \).

Secondly, in the class of random dynamical systems generated by random differential equations, we obtained both global and local
\( C^m \)-linearization results. These results provide the foundation to extend our analysis to the class of random dynamical systems generated by stochastic differential equations (SDEs). To achieve smooth linearization results for random dynamical systems generated by SDEs, we leveraged the connection between random differential equations (RDEs) and SDEs established in \cite{imkeller}.

The  paper is arranged as follows. In the remaining part of this section, we introduce the basic notations and definitions needed throughout the paper. Section 2 focuses on semilinear Carathéodory equations, where we establish the smooth topological equivalence results, including sufficient conditions for \( C^m \)-linearizations and detailed properties of the conjugating map. In Section 3, we apply these results to the setting of random dynamical systems. Specifically, we analyze systems generated by random differential equations (RDEs) and stochastic differential equations (SDEs), leveraging the connection between RDEs and SDEs to derive smooth \( C^m \)-linearization results for the corresponding random dynamical systems. The Appendix contains some key theorems used in proving our results.

\subsection{Preliminaries}
Throughout the paper, let $I\subseteq\mathbb{R}$ be an interval unbounded from above, $X$ be a Banach space equipped with a norm $\norm{\cdot},$ $L(X)$ is the space of linear operator on $X,$ and $(\Omega,\fF)$ be a measurable space. Moreover, let $\{\norm{\cdot}_{t,\omega}\}_{t\in I,\,\omega\in\Omega}$ denote a family of norms on $X$ such that
\begin{enumerate}
	\item[(i)] the mapping $(t,\omega,x)\mapsto\norm{x}_{t,\omega}$ is strongly measurable,

	\item[(ii)] there exists a measurable mapping $\ell:I\times\Omega\to(0,\infty)$ such that $\ell(\cdot,\omega)$ is locally bounded  and
	\begin{equation}
		\frac{1}{\ell(t,\omega)}\norm{x}\leq\norm{x}_{t,\omega}\leq\ell(t,\omega)\norm{x}
		\fall t\in I,\,\omega\in\Omega,\,x\in X.
		\label{normest}
	\end{equation}
\end{enumerate}
Furthermore, for a linear operator $T\in L(X)$ we  define
$$
	\norm{T}_{s,t,\omega}:=\sup_{\norm{x}_{s,\omega}=1}\norm{Tx}_{t,\omega}\fall s,t\in I,\,\omega\in\Omega,
$$
and obtain $\norm{Tx}_{t,\omega}\leq\norm{T}_{s,t,\omega}\norm{x}_{s,\omega}$ for all $t,s\in I$ and $x\in X$.

Given this, a function $\mu:I\to X$ is called  \emph{bounded} (for $\omega\in\Omega$), if
$$
	\norm{\mu}_\omega:=\sup_{t\in I}\norm{\mu(t)}_{t,\omega}<\infty.
$$
We write $BC_\omega(I,X)$ for the set of all bounded and continuous functions $\mu:I\to X$, which equipped with the norm $\norm{\cdot}_\omega$ becomes a Banach space.

We denote $f:I\times X\to Y$ as \emph{Carath{\'e}odory function}, $Y$ being another Banach space, provided the following hold:
\begin{itemize}
	\item $f(\cdot,x):I\to Y$ is strongly measurable for each $x\in X$ (w.r.t.\ the Borel $\sigma$-algebras on $I$ and $X$),

	\item $f(t,\cdot):X\to Y$ is continuous for almost all $t\in I.$
\end{itemize}

 Let $f:I\times X\times \Omega\to X$ be a mapping such that for every choice of $\omega\in\Omega$ the mapping $f(\cdot,\cdot,\omega):I\times X\to X$ has the Carath{\'e}odory property. Then the equation
 \begin{equation}\label{par_eq}
   \dot{x} = f(t,x,\omega)
 \end{equation}
 is called an \emph{ $\omega-$ dependent Carath{\'e}odory differential equation,} and an absolutely continuous mapping $\phi:J\to X$ is called a \emph{solution} of \eqref{par_eq} (to the parameter $\omega$), if the mapping $f(\cdot,\phi(\cdot),\omega)$ is locally integrable and the identity
 \begin{equation}\label{sol_par}
  \phi(t) = \phi(\tau_0) + \int_{\tau_0}^{t} f(s,\phi(s),\omega)ds
 \end{equation}
 holds for all $\tau_0, t\in J,$ where $J\subseteq I$ is a subinterval. If in addition the mapping $\phi$ satisfies the condition $\phi(\tau) = \xi$ then we say that $\phi$ solves the initial value problem for $\tau,\xi\in X$.

 The next lemma bridges the integral and differential formulations of solutions to the Carath{\'e}odory equation, providing a tool for analyzing these solutions within a functional-analytic framework.

\begin{lemma}
    An absolutely continuous mapping $\phi:J\to X$ solves \eqref{par_eq} if and only if $\dot{\phi}(t) = f(t,\phi(t),\omega)$ for almost all $t\in J.$
\end{lemma}

 The following theorem presents the results on existence, uniqueness and continuous dependence of solutions on initial values and parameters.
\begin{theorem}\label{eq-uni}\cite[Thm. 2.4, p. 49]{Wanner_96}, \cite{bres_07}
Let $f:I\times X\times \Omega\to X$ be a mapping such that, for every  $\omega\in\Omega,$ the mapping $f(\cdot,\cdot,\omega):I\times X\to X$ satisfies the Carath{\'e}odory property. Furthermore, let $l, l_0:I\to\mathbb{R}_0^+$ be locally integrable functions such that the two estimates
\begin{equation*}
  \|f(t,x,\omega)-f(t,\bar{x},\omega)\|\leq l(t)\|x-\bar{x}\|,
\end{equation*}
\begin{equation*}\label{bound_par}
  \|f(t,0,\omega)\|\leq l_0(t)
\end{equation*}
hold for almost all $t\in I$ and for all $x,\bar{x}\in X$ and $\omega\in\Omega.$

Then for any $(\tau,\xi,\omega_0)\in I\times X\times \Omega$ the differential equation \eqref{par_eq} has a unique solution $\phi(\cdot;\tau,\xi,\omega_0):I\to X$ to the parameter value $\omega_0,$ satisfying the initial condition $x(\tau) = \xi,$ and the so-defined mapping $\phi(\cdot,\cdot,\cdot,\omega):I\times I\times X\to X$ is continuous.
\end{theorem}
The mapping $\varphi:I\times I\times X\times \Omega\to X$ provided by Theorem \ref{eq-uni} is called a \emph{general solution} of differential equation \eqref{par_eq}. For all $\tau, t_1, t_2\in I, \xi\in X$ and $\omega\in\Omega,$ this mapping satisfies the two-parameter group property
\begin{equation}\label{cocycle_prop_par}
  \varphi(t_2;\tau,\xi,\omega) = \varphi(t_2;t_1,\varphi(t_1;\tau,\xi,\omega),\omega).
\end{equation}

\section{Semilinear  Carath{\'e}odory equations}
In this section, we focus on semilinear Carath{\'e}odory differential equations with  parameter  and investigate their topological equivalence to linear systems. Specifically, we aim to derive sufficient conditions for smooth linearization of class \( C^m \).

Let $I\subseteq\mathbb{R}$ be an interval unbounded from above and $(\Omega,\fF)$ be a measurable space, $A_{\omega}:I\to L(X)$ be locally integrable  and $F_{\omega}:I\times X\to X$ with a given $\omega\in\Omega.$

Consider the system of semilinear  Carath{\'e}odory differential equations with  parameter $\omega\in\Omega$
\begin{equation}\tag{$N$}
	\dot x=A_\omega(t)x+F_\omega(t,x).
	\label{deq}
\end{equation}
Together with the system \eqref{deq} consider a linear system
\begin{equation}\tag{$L$}
    \dot{x} = A_{\omega}(t)x.
    \label{lin}
\end{equation}
\begin{lemma}\cite[Lem. 2.9, p. 56]{Wanner_96}\label{evol} Let $\phi_{\omega}:I\times I\times X\to X$ denote the general solution of the linear equation \eqref{lin}. Then for any $s,t\in I$ the mapping $\Phi_{\omega}:I\times I\to L(X),$
\begin{equation*}
  \Phi_{\omega}(t,s):=\phi_{\omega}(t;s,\cdot),
\end{equation*}
 is called the evolution operator of \eqref{lin} and satisfies the following properties:
\begin{itemize}
  \item[(i)] each $\Phi_{\omega}(t,s),$ $s,t\in I$ is an element of the group $GL(X)$ of invertible operators in $L(X)$ and the mapping $\Phi_{\omega}:I\times I\to GL(X)$ is continuous;
  \item[(ii)] the identities
  \begin{align*}
    \Phi_{\omega}(t,t) =Id_{X}, \\
    \Phi_{\omega}(t,s) = \Phi_{\omega}(t,r)\Phi_{\omega}(r,s),\\
    \Phi_{\omega}(t,s)^{-1} = \Phi_{\omega}(s,t)
  \end{align*}
  hold for all $r,s,t\in I.$
\end{itemize}
\end{lemma}
\subsection{Topological linearization}
Assume that the following assumptions hold for all $\omega\in\Omega$
\begin{enumerate}
	\item[$(H_1)$] $A_\omega:I\to L(X)$ is locally integrable and there exist reals $K(\omega)\geq 1$ and $\alpha>0$ such that	
	\begin{equation}
		\norm{\Phi_\omega(t,s)}_{s,t,\omega}\leq K(\omega)e^{-\alpha(t-s)}\fall t,s\in I,\,s\leq t.
		\label{hyp11}
	\end{equation}

	\item[$(H_2^0)$] $F_\omega:I\times X\to X$ is a Carath{\'e}odory function and there exist reals $L(\omega), M(\omega)\geq 0$ such that
	\begin{align}
		F_\omega(t,0)&=0,
		\label{hyp20}\\
		\norm{F_\omega(t,x)}_{t,\omega}&\leq M(\omega), \label{hyp21}\\
  \norm{F_\omega(t,x)-F_\omega(t,\bar x)}_{t,\omega}&\leq L(\omega)\norm{x-\bar x}_{t,\omega}
		\label{hyp22}
	\end{align}
 for almost all $t\in I$ and all $x,\bar x\in X.$
\end{enumerate}

It is easy to verify that under assumptions $(H_1),$ $(H_2^0),$ the semilinear  equation \eqref{deq} satisfies the assumptions of Theorem~\ref{eq-uni} with $f(t,x,\omega) = A_\omega(t)x+F_\omega(t,x),$ $l(t):= \|A_{\omega}(t)\|+L(\omega)$ and $l_0(t) = 0.$  Then, for every initial condition $x(\tau) = \xi$ and fixed $\omega\in\Omega$ there exists exactly one solution of \eqref{deq} on $I$. From now, we denote the general solution of \eqref{deq} by $\varphi_{\omega}.$

Before delving into the main result, let us first introduce some technical lemmas. Throughout the paper, $\tau_0\in\mathbb{R}$ is fixed and $\tau\in I:=[\tau_0,\infty)$.
\begin{lemma}[the operator $\cL_{\omega}$]
	 If $(H_1)$ holds, then
    \begin{align*}
        \cL_{\omega}:I\times X &\to BC_\omega(I,X), &
        \cL_{\omega}(\tau,\xi) &:= \Phi_\omega(\cdot,\tau)\xi
    \end{align*}
   defines an operator such that $L_\omega(\tau,\cdot)$ is bounded and linear for all $\tau\in I$, $\omega\in\Omega.$ Moreover, for all $\omega \in \Omega$ and $\xi \in X$, the following holds:
    \begin{itemize}
        \item[(i)] $(\tau, \xi) \mapsto \cL_{\omega}(\tau,\xi)$ is continuous,
        \item[(ii)] $\cL_{\omega}(\tau,\xi):I \to X$ is absolutely continuous.
    \end{itemize}
\end{lemma}\label{leml}
\begin{proof}
	Let $\omega\in\Omega$, $\tau\in I,$ and $\xi\in X$. Then $\Phi_\omega(\cdot,\tau)\xi:I\to X$ solves the initial value problem $\dot x=A_\omega(t)x$, $x(\tau)=\xi,$ and is therefore absolutely continuous.

   Without loss of generality, let $\tilde{\tau} \leq \tau \leq t$. To prove the continuity of the map $(\tau, \xi) \mapsto \cL_{\omega}(\tau,\xi)$, consider the following expression:
\begin{align*}
    \norm{\cL_{\omega}(\tau,\xi) - \cL_{,\omega}(\tilde{\tau},\xi_0)}_{t,\omega}
    &\leq \norm{\Phi_\omega(t,\tau)\xi - \Phi_\omega(t,\tau)\xi_0}_{t,\omega}
    + \norm{\Phi_\omega(t,\tau)\xi_0 - \Phi_\omega(t,\tilde{\tau})\xi_0}_{t,\omega} \\
    &\leq K(\omega)e^{-\alpha(t-\tau)}\norm{\xi-\xi_0}
    + \norm{\Phi_\omega(t,\tau) - \Phi_\omega(t,\tilde{\tau})}_{t,\omega}\norm{\xi_0}.
\end{align*}

From the continuous dependence of $\Phi_\omega(t,\tau)$ on the initial time $\tau$, it follows that for any $\varepsilon > 0$, there exists $\delta > 0$ such that if $|\tau -\tilde{\tau}| + \norm{\xi - \xi_0} < \delta$, then
\[
\norm{\cL_{\omega}(\tau,\xi) - \cL_{\omega}(\tilde{\tau},\xi_0)}_{t,\omega}
\leq K(\omega)\delta + \norm{\xi_0}\varepsilon \leq \tilde{\varepsilon},
\]
where $\tilde{\varepsilon}$ can be made arbitrarily small. Therefore, the continuity of $(\tau, \xi) \mapsto \cL_{\omega}(\tau,\xi)$ follows.

    Moreover, one has
	\begin{align}\label{norm_cL}
		\norm{[\cL_{\omega}(\tau,\xi)](t)}_{t,\omega}
		&=
		\norm{\Phi_\omega(t,\tau)\xi}_{t,\omega}
		\stackrel{\eqref{hyp11}}{\leq}
		K(\omega)e^{-\alpha(t-\tau)}\norm{\xi}_{\tau,\omega}
		=
		K(\omega)e^{-\alpha(t-\tau_0)}e^{\alpha(\tau-\tau_0)}\norm{\xi}_{\tau,\omega}\nonumber\\
		&\leq
		K(\omega)e^{\alpha(\tau-\tau_0)}\norm{\xi}_{\tau,\omega}
		\stackrel{\eqref{normest}}{\leq}
		K(\omega)e^{\alpha(\tau-\tau_0)}\ell(\tau,\omega)\norm{\xi}\fall\tau_0\leq t
	\end{align}
	and passing to the supremum over $t\in I$ yields $\norm{\cL_{\omega}(\tau,\xi)}_\omega\leq K(\omega)e^{\alpha(\tau-\tau_0)}\ell(\tau,\omega)\norm{\xi}$, which is the claim.
\end{proof}

\begin{lemma}[the operator $\cF_\omega$]
	\label{lemf}
	If $(H_1)$ and $(H_2^0)$ hold, then the nonlinear operator
	\begin{align}
		\cF_\omega:BC_\omega(I,X)&\to BC_\omega(I,X),&
		\cF_\omega(\phi)&:=\int_{\tau_0}^\cdot\Phi_\omega(\cdot,s)F_\omega(s,\phi(s))\, ds
	\end{align}
	is well-defined,  and, for all $\phi,\bar\phi\in BC_\omega(I,X),\,\omega\in\Omega$ satisfies
	\begin{equation}\label{eq2.8}
		\norm{\cF_\omega(\phi)-\cF_\omega(\bar\phi)}_\omega
		\leq
		\frac{K(\omega)L(\omega)}{\alpha}\norm{\phi-\bar\phi}_\omega,
	\end{equation}
 \begin{equation}\label{eq2.9}
		\norm{\cF_\omega(\phi)}_\omega
		\leq
		\frac{K(\omega)M(\omega)}{\alpha}.
	\end{equation}
\end{lemma}
\begin{proof}
	Let $\omega\in\Omega$ and $\phi,\bar\phi\in BC_\omega(I,X).$ Then $\displaystyle\int_{\tau_0}^{\cdot}\Phi_\omega(\cdot,s)F_\omega(s,\phi(s))\, ds$ solves the linearly inhomogeneous initial value problem $\dot x=A_\omega(t)x + F_\omega(t,\phi(t))$, $x(\tau_0)=0.$   Moreover,
\begin{align*}
		\norm{\cF_\omega(\phi)(t)}_{t,\omega}
		&{\leq} \int_{\tau_0}^t\norm{\Phi_\omega(t,\tau)}_{s,t,\omega}\norm{F_{\omega}(s,\phi(s))}_{s,\omega}\, ds
        \stackrel{\eqref{hyp11}}{\leq} K(\omega) \int_{\tau_0}^te^{-\alpha(t-s)}\norm{F_{\omega}(s,\phi(s))}_{s,\omega}\, ds
        \\& \stackrel{\eqref{hyp21}}{\leq} K(\omega) M(\omega) \int_{\tau_0}^te^{-\alpha(t-s)}\, ds \leq \frac{K(\omega)M(\omega) }{\alpha}\fall\tau_0\leq t,
	\end{align*}
and passing to the supremum over $t\in I$ yields $\cF_\omega(\phi)\in BC_\omega(I,X).$

Furthermore,
\begin{align*}
		\norm{\cF_\omega(\phi)(t)-\cF_\omega(\bar\phi)(t)}_{t,\omega}
		&{\leq} \int_{\tau_0}^t\norm{\Phi_\omega(t,\tau)}_{s,t,\omega}\norm{F_{\omega}(s,\phi(s))- F_{\omega}(s,\bar\phi(s))}_{s,\omega}\, ds\\
        &{\stackrel{\eqref{hyp11}}{\leq}} K(\omega) \int_{\tau_0}^te^{-\alpha(t-s)}\norm{F_{\omega}(s,\phi(s))- F_{\omega}(s,\bar\phi(s))}_{s,\omega}\, ds
        \\ &{\stackrel{\eqref{hyp22}}{\leq}}K(\omega) L(\omega) \int_{\tau_0}^te^{-\alpha(t-s)}\norm{\phi(s)-\bar\phi(s)}_{s,\omega}\, ds
        \\ &{\leq} \frac{K(\omega)L(\omega)}{\alpha}\norm{\phi-\bar\phi}_{\omega}\fall\tau_0\leq t.
	\end{align*}
Passing to the supremum over $t\in I$ yields \eqref{eq2.8} for all \mbox{$\phi,\bar\phi\in BC_\omega(I,X)$}.
Finally, $\cF(\phi)$ is absolutely continuous.
\end{proof}

\begin{theorem}[topological linearization of \eqref{deq}]\label{thm_lin}

	If $(H_1)$ and $(H_2^0)$ hold with
	$$
		\displaystyle
K(\omega)L(\omega)<\alpha,$$
	then there exists a   map $H_{\omega}:I\times X\to X$ such that, for every $\omega\in\Omega,$ $s,t\in I$:
\begin{itemize}
  \item[(i)] $H_\omega(t,\cdot)$ is a homeomorphism on $X$ with $H_{\omega}(t,0) = 0;$
  \item[(ii)] the linear equation \eqref{lin} and the semilinear equation \eqref{deq} are conjugated in the sense that
  \begin{equation*}
    H_\omega(t,\Phi_{\omega}(t,s)\cdot) =\varphi_{\omega}(t,s, H_{\omega}(s,\cdot));
  \end{equation*}
   \item[(iii)] the operators $H_{\omega}(t,\cdot)$ and its inverse are near identity, i.e. for all $\xi\in X,$ $\eta\in X$:
   $$\norm{H_{\omega}(t,\xi)-\xi}_{t,\omega}\leq \frac{K(\omega)M(\omega)}{\alpha} \quad\text{and}\quad \norm{H_{\omega}(t,\cdot)^{-1}(\eta)-\eta}_{t,\omega}\leq \frac{K(\omega)M(\omega)}{\alpha};$$
   \item[(iv)]  the operators $H_{\omega}(t,\cdot):X\to X$ and its inverse satisfy a global Lipschitz condition;
  \item[(v)] the operators $H_{\omega}:I\times X\to X$ and $(t,x)\mapsto H_{\omega}^{-1}(t,\cdot)(x)$  are continuous.
\end{itemize}
We denote $H_\omega$ as a (topological) conjugacy between the semilinear equation \eqref{deq} and its linear part \eqref{lin}.
\end{theorem}
\begin{proof}
For clarity of presentation, we carry out the proof in four steps.

\emph{Step 1.}

	We fix $\tau\in I, \omega\in \Omega.$ Let us define  $\cT_{\omega}:BC_\omega(I,X)\tm X\to BC_\omega(I,X)$ as
	\begin{equation}\label{oper_tau}
		\cT_{\omega}(\phi;\xi,\tau):=\cF_\omega(\phi+\cL_{\omega}(\tau,\xi)),
	\end{equation}
where $\cF_\omega:BC_\omega(I,X)\to BC_\omega(I,X)$ is introduced in Lemma \ref{lemf} and $\cL_{\omega}:X\to BC_\omega(I,X)$ is introduced in Lemma \ref{leml}. We have

\begin{align*}
  \norm{\cT_{\omega}(\phi;\xi,\tau)-\cT_{\omega}(\bar{\phi};\xi,\tau)}_{\omega}&\overset{\eqref{oper_tau}}{=}\norm{\cF_\omega(\phi+\cL_{\omega}(\tau,\xi))-\cF_\omega(\bar{\phi}+\cL_{\omega}(\tau,\xi))}_{\omega}\\
		&{\stackrel{\eqref{eq2.8}}{\leq}}\frac{K(\omega)L(\omega)}{\alpha}\norm{\phi-\bar{\phi}}_{\omega}
\end{align*}
for all \mbox{$\phi,\bar\phi\in BC_\omega(I,X).$}
 Hence, $\cT_{\omega}$ is a contraction uniformly in $\xi\in X$ and in $\tau\in I$ and by uniform contraction principle (Theorem~\ref{uniform contraction}) there is a unique fixed point $\varphi^*_{\tau,\omega}(\xi)\in BC_{\omega}(I,X)$ of $\cT_{\omega}(\cdot;\xi,\tau)$ for all $\xi\in X,$ $\tau\in I$ and the mapping $(\tau,\xi)\mapsto\varphi^*_{\tau,\omega}(\xi)$ is continuous.

 Moreover,
 \begin{align}
  \norm{\cT_{\omega}(\phi;\xi,\tau)-\cT_{\omega}(\phi;\bar{\xi},\tau)}_{\tau,\omega}&\overset{\eqref{oper_tau}}{=}\norm{\cF_\omega(\phi+\cL_{\omega}(\tau,\xi))-\cF_\omega(\phi+\cL_{\omega}(\tau,\bar{\xi}))}_{\tau,\omega}\nonumber\\
		&{\stackrel{\eqref{eq2.8}}{\leq}}\frac{K(\omega)L(\omega)}{\alpha}\norm{\cL_{\omega}(\tau,\xi)-\cL_{\omega}(\tau,\bar{\xi})}_{\tau,\omega}\nonumber\\&{\stackrel{\eqref{norm_cL}}{\leq}}\frac{K^2(\omega)L(\omega)}{\alpha}e^{\alpha(\tau-\tau_0)}\ell(\tau,\omega)\norm{\xi-\bar{\xi}}.
\end{align}

This leads us to conclude that $\cT_{\omega}(\phi,\cdot,\tau)$ is Lipschitz with Lipschitz constant $\displaystyle\frac{K^2(\omega)L(\omega)}{\alpha}e^{\alpha(\tau-\tau_0)}\ell(\tau,\omega).$ Therefore,  the mapping $\xi\mapsto\varphi^*_{\tau,\omega}(\xi)$ is also  Lipschitz with Lipschitz constant $\displaystyle\frac{K^2(\omega)L(\omega)e^{\alpha(\tau-\tau_0)}\ell(\tau,\omega)}{\alpha-{K(\omega)L(\omega)}}.$

Furthermore, we show that $\varphi^*_{\tau,\omega}(\xi)$ satisfies   the identity
\begin{equation}
    \varphi^*_{\tau,\omega}(\xi) = \varphi^*_{r,\omega}(\Phi_{\omega}(r,\tau)\xi)
\end{equation}
 for all $\xi\in X$ and almost all $\tau, r\in I.$

Indeed,
$$\varphi^*_{\tau,\omega}(\xi) = \int_{\tau_0}^t\Phi_{\omega}(t,s)F_{\omega}(s,\varphi_{\tau,\omega}^*(\xi)(s)+\Phi_{\omega}(t,\tau)\xi)\,ds$$
is a  unique solution to the initial value problem $\dot x = A_{\omega}(t)x + F_{\omega}(t, x + \Phi_{\omega}(t,\tau) \xi),
x(\tau_0) = 0.
$
On the other hand,
\begin{eqnarray*}
    \varphi^*_{r,\omega}(\Phi_{\omega}(r,\tau)\xi) &=& \int_{\tau_0}^t\Phi_{\omega}(t,s)F_{\omega}(s,\varphi^*_{r,\omega}(\Phi_{\omega}(r,\tau)\xi)(s)+\Phi_{\omega}(t,r)\Phi_{\omega}(r,\tau)\xi)\,ds\\
    &=& \int_{\tau_0}^t\Phi_{\omega}(t,s)F_{\omega}(s,\varphi^*_{r,\omega}(\Phi_{\omega}(r,\tau)\xi)(s)+\Phi_{\omega}(t,\tau)\xi)\,ds
\end{eqnarray*}
is a  unique solution to the same initial value problem. As consequence of uniqueness of solution we obtain the desired result.

\emph{Step 2.}
 Let us introduce  $H_{\omega}:I\times X\to X$ and  $G_{\omega}:I\times X\to X$ as follows
 \begin{align}\label{H}H_\omega{(t,\xi)}&{:=}\xi+\int_{\tau_0}^t \Phi_{\omega}(t,s)F_{\omega}(s,\varphi^*_{t,\omega}(\xi)(s)+\Phi_{\omega}(t,s)\xi)\,ds\nonumber\\
 &{=}\xi+\varphi^*_{t,\omega}(\xi)(t).\\
 G_\omega{(t,\eta)}&{:=}\eta - \int_{\tau_0}^t \Phi_{\omega}(t,s)F_{\omega}(s,\varphi_{\omega}(s,t,\eta))\,ds.\nonumber
 \end{align}

 We can show the boundedness of $\xi\mapsto H_\omega{(t,\xi)} - \xi$ and $\xi\mapsto G_\omega{(t,\eta)} -\eta$ as follows
 \begin{align*}
     \norm{H_\omega{(t,\xi)} - \xi}_{t,\omega}&{\leq}\int_{\tau_0}^t\norm{\Phi_{\omega}(t,s)}_{s,t,\omega}\norm{F_{\omega}(s,\varphi^*_{t,\omega}(\xi)(s)+\Phi_{\omega}(t,s)\xi)}_{s,\omega} \,ds\\
     &{\stackrel{\eqref{hyp11}}{\leq}}K(\omega)\int_{\tau_0}^t e^{-\alpha(t-s)}\norm{F_{\omega}(s,\varphi^*_{t,\omega}(\xi)(s)+\Phi_{\omega}(t,s)\xi)}_{s,\omega} \,ds\\
     &{\stackrel{\eqref{hyp21}}{\leq}} K(\omega)M(\omega)\int_{\tau_0}^t e^{-\alpha(t-s)}\,ds \leq \frac{K(\omega)M(\omega)}{\alpha}\fall\xi\in X,
 \end{align*}
 and respectively,
 \begin{align*}
    \norm{G_\omega{(t,\eta)} - \eta}_{t,\omega}&{\leq}\int_{\tau_0}^t\norm{\Phi_{\omega}(t,s)}_{s,t,\omega}\norm{F_{\omega}(s,\varphi_{\omega}(s,t,\eta)}_{s,\omega} ds\\
     &{\leq}\frac{K(\omega)M(\omega)}{\alpha}\fall\eta\in X.
 \end{align*}

Now, let us demonstrate the conjugation of \eqref{deq} and \eqref{lin}. First, we claim that \( H_\omega(t,\Phi_{\omega}(t,\tau)\xi) \) is differentiable in \( t \). Indeed,
\[
H_\omega(t,\Phi_{\omega}(t,\tau)\xi) = \Phi_{\omega}(t,\tau)\xi + \int_{\tau_0}^t \Phi_\omega(t,s) F_\omega\left(s, \varphi^*_{t,\omega}(\Phi_\omega(t,\tau)\xi)(s)
+ \Phi_\omega(t,s)\Phi_\omega(t,\tau)\xi \right) \, ds.
\]

The function \( \Phi_\omega(t,s) \) is the solution of \eqref{lin}, and the conditions in \( (H_1) \) guarantee that \( \Phi_\omega \) is continuously differentiable with respect to \( t \) for \( t \geq s \). Specifically:
\[
\frac{d}{d t} \Phi_\omega(t,s) = A_\omega(t)\Phi_\omega(t,s).
\]

To analyze the second term, we note that \( (H_2^0) \) guarantees that \( F_\omega \) is sufficiently regular for the integral term to be differentiable. Moreover, \( \varphi^*_{t,\omega}(\xi) \) solves \eqref{deq} and depends on \( t \) in a smooth way through \( \Phi_\omega(t,s) \) and the nonlinearity \( F_\omega \). The Lipschitz continuity of \( F_\omega \) ensures that the solution \( \varphi^*_{t,\omega}(\xi) \) is differentiable with respect to \( t \).

Combining all this with the differentiation under the integral Theorem \ref{diff_int}, we can compute the derivative in a direct form:
\begin{eqnarray*}
        \frac{d }{d t}H_\omega(t,\Phi_{\omega}(t,\tau)\xi)
        &=&\frac{d }{d t}\Phi_{\omega}(t,\tau)\xi + \frac{d }{d t}\varphi^*_{t,\omega}(\Phi_{\omega}(t,\tau)\xi)(t)\\
        &=&A_\omega(t)\left(\Phi_{\omega}(t,\tau)\xi +  \varphi^*_{t,\omega}(\Phi_{\omega}(t,\tau)\xi)(t)\right) + F_\omega \left(t,H_\omega(t,\Phi_{\omega}(t,\tau)\xi)\right) \\
        &=& A_\omega(t)H_\omega(t,\Phi_{\omega}(t,\tau)\xi) + F_\omega\left(t, H_\omega(t,\Phi_{\omega}(t,\tau)\xi)\right).
\end{eqnarray*}

This means that \( t \mapsto H_{\omega}(t,\Phi_{\omega}(t,\tau)\xi) \) solves the initial value problem:
\[
\dot{x} = A_{\omega}(t)x + F_{\omega}(t,x), \quad x(\tau) = H_\omega (\tau,\xi).
\]

By the uniqueness of the solution (due to Theorem \ref{eq-uni}), we obtain:
\[
H_{\omega}(t,\Phi_{\omega}(t,\tau)\xi) = \varphi_{\omega}(t,\tau,H_{\omega}(\tau,\xi)) \quad \forall \, \xi \in X.
\]
 In the same way it can be shown that
\begin{equation*}
  G_{\omega}(t,\varphi_{\omega}(t,s,\eta)) = \Phi_{\omega}(t,s)G_{\omega}(s,\eta) \fall\eta\in X.
\end{equation*}

\emph{Step 3.} In order to complete the proof and to show that $H_{\omega}(t,\cdot)$ is a homeomorphism, we start with the claim that $H_{\omega}(t,\cdot):X\to X$ is bijective. Indeed, for any $\eta\in X,$ we have

\begin{align*}
     H_\omega(t,G_{\omega}(t,\varphi_{\omega}(s,t,\eta))) &\overset{\eqref{H}}{=}G_{\omega}(t,\varphi_{\omega}(t,\tau,\eta))
     +\int_{\tau_0}^t\Phi_{\omega}(t,s)F_{\omega}(s,\varphi^*_{t,\omega}(G_{\omega}(s,\varphi_{\omega}(t,\tau,\eta)))(s)\\
     &{+}\Phi_{\omega}(t,s)G_{\omega}(s,\varphi_{\omega}(t,\tau,\eta))))\, ds\\
     &{=}\varphi_{\omega}(t,\tau,\eta)-\int_{\tau_0}^t \Phi_{\omega}(t,s)F_{\omega}(s,\varphi_{\omega}(t,\tau,\eta))ds\\
     &{+}\int_{\tau_0}^t\Phi_{\omega}(t,s)F_{\omega}(s,\varphi^*_{t,\omega}(G_{\omega}(s,\varphi_{\omega}(t,\tau,\eta)))(s)\\
     &{+}
     \Phi_{\omega}(t,s)G_{\omega}(s,\varphi_{\omega}(t,\tau,\eta))))\, ds.
 \end{align*}
 If we abbreviate by $\nu(t) = \norm{H_\omega(t,G_{\omega}(t,\varphi_{\omega}(s,t,\eta)))  - \varphi_{\omega}(t,\tau,\eta)}_{t,\omega},$ then we have
 \begin{align*}
     \nu(t) &{\leq}\int_{\tau_0}^t \norm{ \Phi_{\omega}(t,s)}_{s,t,\omega}\|F_{\omega}(s,\varphi^*_{t,\omega}(G_{\omega}(s,\varphi_{\omega}(t,\tau,\eta)))(s)\\
     &{+}\Phi_{\omega}(t,s)G_{\omega}(s,\varphi_{\omega}(t,\tau,\eta)))- F_{\omega}(s,\varphi_{\omega}(t,\tau,\eta))\|_{s,\omega}ds\\
     &{\leq}K(\omega)L(\omega)\int_{\tau_0}^{t}e^{-\alpha(t-s)}\norm{\varphi^*_{t,\omega}(G_{\omega}(s,\varphi_{\omega}(t,\tau,\eta)))(s)+\Phi_{\omega}(t,s)G_{\omega}(s,\varphi_{\omega}(t,\tau,\eta)))- \varphi_{\omega}(t,\tau,\eta)}_{s,\omega}\, ds\\
     &{\leq}K(\omega)L(\omega)\int_{\tau_0}^t e^{-\alpha(t-s)}\norm{H_{\omega}(s,\Phi_{\omega}(t,s)G_{\omega}(s,\varphi_{\omega}(t,\tau,\eta))))-\varphi_{\omega}(t,\tau,\eta)}_{s,\omega}\, ds\\
     &{\leq} K(\omega)L(\omega)\int_{\tau_0}^t e^{-\alpha(t-s)}\nu(s)\, ds \leq \frac{K(\omega)L(\omega)}{\alpha}\sup_{s\in I} \nu(s) \fall t\in I.
     \end{align*}

Hence, from the above estimate and the assumption  $\displaystyle
K(\omega)L(\omega)<\alpha,$ it follows that $\nu(t) = 0$ for any $t\in I$ and with $t=\tau,$ $H_{\omega}(\tau,G_{\omega}(\tau,\eta)) = \eta.$

Now we prove that for all $\xi\in X,$ it holds that $G_{\omega
}(t,H_{\omega}(t,\xi)) = \xi.$ We have

\begin{align*}
     G_{\omega}(t,H_{\omega}(t,\Phi_{\omega}(t,\tau)\xi)) &{=} H_{\omega}(t,\Phi_{\omega}(t,\tau)\xi)-\int_{\tau_0}^t\Phi_{\omega}(t,s)F_{\omega}(s,\varphi_{\omega}(s,t,H_{\omega}(t,\Phi_{\omega}(t,\tau)\xi)))\, ds\\
     &{=}\Phi_{\omega}(t,\tau)\xi + \int_{\tau_0}^t\left(F_{\omega}(s,H_{\omega}(s,\Phi_{\omega}(s,\tau)\xi))-F_{\omega}(s,\varphi_{\omega}(s,\tau,H_{\omega}(\tau,\xi)))\right)\, ds\\
     &{=}\Phi_{\omega}(t,\tau)\xi.
 \end{align*}
Now, if we set $t=\tau,$ then $G_{\omega}(\tau, H_{\omega}(\tau,\xi)) = \xi$ holds.

As a consequence, for any $t\in I$, $H_\omega(t,\cdot): X \to X$ is a bijection, and $G_{\omega}(t,\cdot): X\to X$ is its inverse.

\emph{Step 4.}
In this step,  we  show that the conjugacy and its inverse satisfy Lipschitz conditions with Lipschitz constants $L_{H}$ and $L_{G}$ respectively, where
$$L_{G} = 1+\frac{K^2(\omega)L(\omega)}{2\alpha-K(\omega)L(\omega)}  ,\quad\quad L_{H} = 1+\frac{K^2(\omega)L(\omega)}{2\alpha}+\frac{K^3(\omega)L^2(\omega)e^{\alpha(\tau-\tau_0)}\ell(\tau,\omega)}{\alpha(\alpha-{K(\omega)L(\omega)})}.$$

In order to show the claim for the operator $G_\omega,$ we start with additional estimates.

Note that the general solution $\varphi_{\omega}(\cdot,\tau,\eta)$ of \eqref{deq}  also solves the inhomogeneous equation $$\dot{x} = A_{\omega}(t)x+F_{\omega}(t,\varphi_{\omega}(t,\tau,\eta)).$$
Hence, by variation of constants formula (Theorem \ref{var_const}), we have
$$\varphi_{\omega}(s,t,\eta) = \Phi_{\omega}(t,s)\eta+\int_{s}^{t}\Phi_{\omega}(t,\tau)F_{\omega}(\tau,\varphi_{\omega}(t,\tau,\eta))\,d\tau,$$
so we obtain
 \begin{align*}
\norm{\varphi_{\omega}(s,t,\eta)-\varphi_{\omega}(s,t,\bar{\eta})}_{t,\omega}&{\leq}\norm{\Phi_{\omega}(t,s)(\eta-\bar{\eta})}_{\omega}\nonumber\\  &+\int_{s}^t\norm{\Phi_{\omega}(t,\tau)}_{s,t,\omega}\norm{F_{\omega}(\tau,\varphi_{\omega}(t,\tau,\eta))-F_{\omega}(\tau,\varphi_{\omega}(t,\tau,\bar{\eta}))}_{s,\omega} ds\nonumber\\&\leq
K(\omega)e^{-\alpha(t-s)}\norm{\eta-\bar{\eta}}\\&+K(\omega)L(\omega) \int_{s}^{t}e^{-\alpha (t-\tau)}
     \norm{\varphi_{\omega}(t,\tau,\eta)-\varphi_{\omega}(t,\tau,\bar{\eta})}_{t,\omega}ds .
 \end{align*}
Multiplying both sides by
 $e^{\alpha(t-s)}$ and using Gronwall's
inequality \eqref{Gron1}, gives us
\begin{align*}
\norm{\varphi_{\omega}(s,t,\eta)-\varphi_{\omega}(s,t,\bar{\eta})}_{t,\omega}e^{\alpha(t-s)}&{\leq} K(\omega) e^{K(\omega)L(\omega)(t-s)}\norm{\eta-\bar{\eta}}.
\end{align*}
From the estimate above we obtain
\begin{align}\label{Gron}
\norm{\varphi_{\omega}(s,t,\eta)-\varphi_{\omega}(s,t,\bar{\eta})}_{t,\omega}&{\leq} K(\omega) e^{(-\alpha+K(\omega)L(\omega))(t-s)}\norm{\eta-\bar{\eta}}.
\end{align}

   Now, using \eqref{Gron}, we have that
       \begin{align}\label{Lip_G}
            \norm{G_\omega(t,\eta)-G_\omega(t,\bar{\eta})}_{t,\omega}&
            \leq\norm{\eta-\bar{\eta}} + K(\omega)L(\omega)\int\limits_{\tau_0}^t e^{-\alpha(t- s)}\norm{\varphi_{\omega}(t, s, \bar{\eta})-\varphi_{\omega}(t, s, \eta)}_{s,\omega}ds\nonumber\\
            &\leq \norm{\eta-\bar{\eta}}\left(1+K^2(\omega)L(\omega)\int\limits_{\tau_0}^t e^{(-2\alpha+K(\omega)L(\omega))(t-s)}ds\right)\nonumber\\
&\leq\left(1+\frac{K^2(\omega)L(\omega)}{2\alpha-K(\omega)L(\omega)}\right)\norm{\eta-\bar{\eta}} \nonumber\\
            &\leq L_{G}\norm{\eta-\bar{\eta}}.
    \end{align}

Now we show that $H_\omega$ satisfy a global Lipschitz condition. We have
 \begin{eqnarray}\label{lip_H}
       \norm{H_\omega(t,\xi)-H_\omega(t,\bar{\xi})}_{t,\omega}  &\le& \norm{\xi-\bar{\xi}}+K(\omega)L(\omega)\int\limits_{\tau_0}^t e^{-\alpha(t-s)}\left\{\norm{\Phi_{\omega}(t,s)\xi-\Phi_{\omega}(t,s)\widetilde{\xi}}_{s,\omega}\right.\nonumber\\&+&\left.\norm{\varphi^*_{t,\omega}(\xi)-\varphi^*_{t,\omega}(\bar{\xi})}_{\omega}\right\}ds\nonumber\\
       &\le& \norm{\xi-\bar{\xi}}\left(1+K^2(\omega)L(\omega)\int\limits_{\tau_0}^t e^{-2\alpha(t-s)}ds\right)\nonumber\\&+&\frac{K(\omega)L(\omega)}{\alpha}\norm{\varphi^*_{t,\omega}(\xi)-\varphi^*_{t,\omega}(\bar{\xi})}_{\omega}\nonumber\\
       &\le& \norm{\xi-\bar{\xi}}\left(1+K^2(\omega)L(\omega)\int\limits_{\tau_0}^t e^{-2\alpha(t-s)}ds\right)\nonumber\\&+&\frac{K^3(\omega)L^2(\omega)e^{\alpha(\tau-\tau_0)}\ell(\tau,\omega)}{\alpha(\alpha-{K(\omega)L(\omega)})}\norm{\xi-\bar{\xi}}\nonumber\\
       &\le& \left(1+\frac{K^2(\omega)L(\omega)}{2\alpha}+\frac{K^3(\omega)L^2(\omega)e^{\alpha(\tau-\tau_0)}\ell(\tau,\omega)}{\alpha(\alpha-{K(\omega)L(\omega)})}\right)\norm{\xi-\bar{\xi}}\nonumber\\&\leq& L_{H}\norm{\xi-\bar{\xi}}.
    \end{eqnarray}

Based on \eqref{Lip_G} and \eqref{lip_H}, we conclude that $H_\omega(t,\cdot)$ and $G_\omega(t,\cdot)$ are Lipschitz, and therefore uniformly continuous for any $t\in I$.

Moreover, it can be shown that the operator $H_{\omega}:I\times X\to X$ and its inverse are continuous.

    Indeed, from the continuity of the solution of \eqref{lin}, it follows that for any $\varepsilon_1>0$ there exists $ \delta_1(t_0,\xi_0,\varepsilon_1)>0:$ $|t-t_0|+\norm{\xi-\xi_0}<\delta_1\Rightarrow\norm{\Phi_{\omega}(t,s)\xi-\Phi_{\omega}(t_0,s)\xi_0}_{s,\omega}<\varepsilon_1.$

Thus, at the first step we proved the continuity of the mapping $(t,\xi)\mapsto\varphi^*_{t,\omega}(\xi)$. This means that,  for any $\varepsilon_2>0,$ there exists $ \delta_2(t_0,\xi_0,\varepsilon_2)>0:$ $|t-t_0|+\norm{\xi-\xi_0}<\delta_2\Rightarrow\norm{\varphi^*_{t,\omega}(\xi)-\varphi^*_{t,\omega}(\xi_0)}_{t,\omega}<\varepsilon_2.$

   Without loss of generality, let $\tau_0\leq t_0\leq t.$  In order to prove the continuity of the operator $H_{\omega}$, we  consider the following expression

\begin{eqnarray*}
    H_{\omega}(t,\xi)-H_{\omega}(t_0,\xi_0)&{=}&\xi-\xi_{0}+\int_{\tau_0}^{t}\Phi_{\omega}(t,s)F_{\omega}(s,\varphi_{t,\omega}^*(\xi)(s)+\Phi_{\omega}(t,s)\xi)\,ds\nonumber\\&{-}&\int_{\tau_0}^{t_0}\Phi_{\omega}(t_0,s)F_{\omega}(s,\varphi_{t_0,\omega}^*(\xi_0)(s)+\Phi_{\omega}(t_0,s)\xi_0)\,ds
    =\xi-\xi_{0}\nonumber\\&{+}&\int_{\tau_0}^{t_0}(\Phi_{\omega}(t,s)-\Phi_{\omega}(t_0,s))F_{\omega}(s,\varphi_{t,\omega}^*(\xi)(s)+\Phi_{\omega}(t,s)\xi)\,ds\nonumber\\
    &{+}&\int_{\tau_0}^{t_0}\Phi_{\omega}(t_0,s)\left(F_{\omega}(s,\varphi_{t,\omega}^*(\xi)(s)+\Phi_{\omega}(t,s)\xi)\right.\nonumber\\&{-}&\left.F_{\omega}(s,\varphi_{t_0,\omega}^*(\xi_0)(s)+\Phi_{\omega}(t_0,s)\xi_0)\,ds\right)\nonumber\\&{+}&\int_{t_0}^{t}
    \Phi_{\omega}(t,s)F_{\omega}(s,\varphi_{t,\omega}^*(\xi)(s)+\Phi_{\omega}(t,s)\xi)\,ds.
\end{eqnarray*}
   Let $C=\max\{\norm{\Phi(t,s)}_{s,t,\omega}:t,s\in I, t\ge s\},$ then we obtain
   \begin{eqnarray*}
    \norm{H_{\omega}(t,\xi)-H_{\omega}(t_0,\xi_0)}_{t,\omega}&{\leq}&\norm{\xi-\xi_{0}}+
M(\omega)\int_{\tau_0}^{t_0}\norm{\Phi_{\omega}(t,s)-\Phi_{\omega}(t_0,s)}_{s,t,\omega}\,ds\nonumber\\
    &{+}&L(\omega)C\int_{\tau_0}^{t_0}\norm{\varphi_{t,\omega}^*(\xi)(s)-\varphi_{t_0,\omega}^*(\xi_0)(s)}_{t,\omega}\,ds\nonumber\\&{+}&M(\omega)C(t-t_0)+L(\omega)C(t_0-\tau_0)\varepsilon\nonumber\\&{\leq}&\norm{\xi-\xi_{0}}+(2M(\omega)+L(\omega)\left(\varepsilon_1+\varepsilon_2)\right)C(t_0-\tau_0)\nonumber\\&+&M(\omega)C(t-t_0),
\end{eqnarray*}
which proves the continuity of $H_{\omega}$ in any $(t_0,\xi_0)\in I\times X.$

In the same way, using the continuity of general solution $\varphi_{\omega}(s,t,\eta)$ of $\eqref{deq},$ we know that for any $\varepsilon_3>0$ there exists $ \delta_3(t_0,\xi_0,\varepsilon_3)>0:$ $|t-t_0|+\norm{\xi-\xi_0}<\delta_3\Rightarrow$ $\norm{\varphi_{\omega}(s,t,\eta)-\varphi_{\omega}(s,t_0,\eta_0)}_{t,\omega}<\varepsilon_3.$ Thus,

\begin{eqnarray*}
    \norm{G_{\omega}(t,\eta)-G_{\omega}(t_0,\eta_0)}_{t,\omega}&{\leq}&\norm{\eta-\eta_{0}}+
M(\omega)\int_{\tau_0}^{t_0}\norm{\Phi_{\omega}(t,s)-\Phi_{\omega}(t_0,s)}_{s,t,\omega}\,ds\nonumber\\
    &{+}&L(\omega)C\int_{\tau_0}^{t_0}\norm{\varphi_{\omega}(s,t,\eta)-\varphi_{\omega}(s,t_0,\eta_0)}_{t,\omega}\,ds+M(\omega)C(t-t_0)\nonumber\\&{\leq}& \norm{\eta-\eta_{0}}+(L(\omega)\varepsilon_3+2M(\omega))C(t_0-\tau_0)\nonumber\\&+&M(\omega)C(t-t_0),
\end{eqnarray*}

which completes the proof.

\end{proof}

  If we set $\Tilde{u}=0$ in \eqref{Gron}, we see that  $\norm{\varphi_{\omega}(t,s,u)}_{t,\omega}\leq K(\omega) e^{(-\alpha+K(\omega)L(\omega))(t-s)}\norm{u}$ which leads us to conclude that the trivial solution of the \eqref{deq} is globally uniformly exponentially stable.

\subsection{Smooth linearization}

We  write $D_2^{j}F_\omega:I\times X\to L_{j}(X)$ for the $j$-th order partial  derivative of a mapping $F_\omega:I\times X\to X$ with respect to the second variable. And we substitute the assumption $(H_2^0)$ by the following condition supposed to hold for some $m\in\N$ and all $\omega\in\Omega$:
\begin{enumerate}
	\item[$(H_2^m)$] $F_\omega:I\times X\to X$ and its partial derivatives $D_2^{j}F_\omega:I\tm X\to L_j(X)$, $1\leq j\leq m$, are Carath{\'e}odory functions and there exist reals $M_j(\omega)\geq 0$ such that
	\begin{align}
		\norm{D_2^{j}F_{\omega}(t,x)}_{t,\omega}&\leq M_j(\omega)\quad\text{for almost all $t\in I$, all $x\in X$ and }0\leq j\leq m.
	\end{align}
\end{enumerate}

\begin{lemma}[Differentiation under the integral]\label{diff}
    Let $m\in\mathbb{N},$ $U\subseteq X, V\subseteq Y$ be open and $f:U\times V \to X$ be such that $f(\cdot,y)$ is integrable for any $y\in V.$ Suppose, that for any $x\in U, y\in V,$ there exist derivatives $D_2^jf(x,y),$ and $$\norm{D_2^jf(x,y)}\leq g(x)$$
    for all $0\leq j\leq m,$ where $g:U\to \mathbb{R}^+$ is integrable.
    Then the mapping $$F(y) = \int_U f(x,y)dx$$ is $m$ times differentiable and its $j$-th derivative is given by$$D^jF(y)= \int_UD^j_2 f(x,y)dx$$
for all $0\leq j\leq m.$
\end{lemma}
%\begin{proof}
%The proof of this statement can be established by mathematical induction.

%First we proof the case $k=1:$
%\begin{eqnarray*}
 %   \lim_{h\to 0}\norm{DF(y) - \int_{U}D_2f(x,y)dx}&=& \lim_{h\to 0}\norm{\frac{F(y+h)-F(y)}{h}-\int_{U}D_2f(x,y)dx}\\&=&\lim_{h\to 0}\norm{\frac{1}{h}\left(\int_{U}f(x,y+h)\,dx - \int_{U}f(x,y)\,dx\right)-\int_{U}D_2f(x,y)dx}\\
  %  &=&\lim_{h\to 0}\norm{\frac{f(x,y+h) - f(x,y)}{h}\,dx - \int_{U}D_2f(x,y)dx}\\
    %&=&\lim_{h\to 0}\norm{\int_{U}D_2f(x,y)dx - \int_{U}D_2f(x,y)dx} = 0.
%\end{eqnarray*}

%Let us assume that for $k=j$ the statement holds.

%Then for $k=j+1$ we have
%\begin{eqnarray*}
 %   D^{j+1}F(y) &=& D(D^jF(y)) = D\left(\int_U D^j_y f(x,y)dx\right) = \int_UD\left(D^j_y f(x,y)\right)dx\\ &=& \int_UD^{j+1}_y f(x,y)dx
%\end{eqnarray*}
%which concludes the proof for any $k\in\mathbb{N}.$
%\end{proof}

\begin{theorem}[$C^m$-linearization of \eqref{deq}]\label{thmdiff}
	Assume  $(H_1)$ and $(H_2^m)$ hold  and  $\displaystyle
K(\omega)M_1(\omega)<\alpha$ for all $\omega\in\Omega.$
	Then all statements of  Theorem \ref{thm_lin} hold with  $H_{\omega}{(t,\cdot)}$ being a $C^m-$ diffeomorphism.
    In the present setting, $H_\omega$ is called a $C^m-$conjugacy between \eqref{deq} and its linearization \eqref{lin}.
\end{theorem}
Note that $C^m$-conjugacies guarantee that the integral manifolds of \eqref{deq} are tangential to the invariant linear integral manifolds of \eqref{lin}.
\begin{proof}

First, we aim to show that
the general solution  $\varphi_{\omega}(s,t,\cdot):X\to X$ of \eqref{deq} is $m$ times differentiable with continuous bounded partial derivatives
$D_3^j\varphi_{\omega}:I\times I\times X\to L_{j}(X)$
for any $\eta\in X,$ $s,t\in I$ and $1\leq j\leq m.$

It is well-known (see, e.g., Theorem 1.16 in \cite{Siegmund}) that  $\varphi_{\omega}\in C^m.$ Our goal is to show the boundedness of the partial derivatives.

From the fact that $\varphi_{\omega}$ is the general solution of \eqref{deq} it follows that, for any $\eta\in X,$ the mapping  $\varphi_\omega$ satisfies the equation $D_1\varphi_{\omega}(t,s,\eta) = A_{\omega}(t)\varphi_{\omega}(t,s,\eta) + F_{\omega}(t,\varphi_{\omega}(t,s,\eta)).$
Then,
\begin{equation*}
  D_1D_3\varphi_{\omega}(t,s,\cdot) = A_{\omega}(t)D_3\varphi_{\omega}(t,s,\eta) + D_2F_{\omega}(t,\varphi_{\omega}(t,s,\eta))D_3\varphi_{\omega}(t,s,\eta)
\end{equation*}
and, therefore, $D_3\varphi_{\omega}(\cdot,s,\eta)$ solves the initial value problem
$\displaystyle\dot{x} = [A_{\omega}(t) +  D_2F_{\omega}(t,\varphi_{\omega}(t,s,\eta))]x,$ $x(s) = Id_{X}.$
Hence, by variation of constants formula we have
\begin{equation*}
  D_3\varphi_{\omega}(t,s,\eta) = \Phi_{\omega}(t,s)+\int_s^t\Phi_{\omega}(r,s) D_2F_{\omega}(r,\varphi_{\omega}(r,s,\eta))D_3\varphi_{\omega}(r,s,\eta)\,dr,
\end{equation*}
and we estimate
\begin{eqnarray*}
    \norm{D_3\varphi_{\omega}(t,s,\eta)}_{t,\omega}&\leq&\int_s^t\norm{\Phi_{\omega}(r,s)}_{s,t,\omega}\norm{D_2F_{\omega}(r,\varphi_{\omega}(r,s,\eta))}_{s,\omega}\norm{D_3\varphi_{\omega}(r,s,\eta)}_{s,\omega}\,dr\\
    &+&\norm{\Phi_{\omega}(t,s)}_{s,t,\omega}\leq K(\omega)e^{-\alpha(t-s)}\\&+&K(\omega)M_1(\omega)\int_s^te^{-\alpha(r-s)}\norm{D_3\varphi_{\omega}(r,s,\eta)}_{s,\omega}\,dr.
\end{eqnarray*}
Multiplying both sides by
 $e^{\alpha(t-s)}$ and using Gronwall's
inequality \eqref{Gron1}, gives us
\begin{equation*}\displaystyle
    \norm{D_3\varphi_{\omega}(t,s,\eta)}_{t,\omega}e^{\alpha(t-s)}\leq K(\omega)e^{K(\omega)M_1(\omega)(t-s)}.
\end{equation*}
From the estimate above, we obtain
\begin{equation*}
||D_3\varphi_{\omega}(t,s,\eta)||_{t,\omega}\leq K(\omega)e^{(K(\omega)M_1(\omega)-\alpha)(t-s)},
\end{equation*}
and the boundedness of the first-order partial derivative follows. Moreover, the boundedness of higher-order partial derivatives can be demonstrated using the same arguments by induction.

Now our goal is to show that the map $\eta\mapsto G_{\omega}{(t,\eta)}$ is  $C^m.$

This follows from the fact that $\varphi_{\omega}(s,t,\cdot):X\to X$ is $C^m,$ hypothesis $(H_2^m),$ and  from applying  Lemma~\ref{diff} to the function $f(t,\eta) = \Phi_{\omega}(t,s)F_{\omega}(s,\varphi_{\omega}(s,t,\eta)).$ To justify the application of Lemma~\ref{diff}, we observe that
\begin{equation*}
    \norm{D_2f(t,\eta)}_{t,\omega}\leq M_1 K^2(\omega) e^{(K(\omega)M_1 - 2\alpha)(t-s)}.
\end{equation*}

This allows us to conclude the desired statement and to express $D_2G_{\omega}(t,\eta)$ in explicit form using the chain rule
\begin{equation}\label{dH}D_2G_{\omega}(t,\eta) = Id_X-\int\limits_{\tau_0}^{t}\Phi_{\omega}(t,s)D_2F_{\omega}(s,\varphi_{\omega}(s,t,\eta))D_3\varphi_{\omega}(s,t,\eta)\,ds.\end{equation}

From $(H_2^m)$ and the fact that $\varphi_{\omega}$ is the general solution of \eqref{deq}, it follows that $D_3\varphi_{\omega}(s,t,\cdot)$ satisfies the differential equation
\begin{equation}\label{dvarphi}
\dot{x} = A_{\omega}(t)x+D_2F_{\omega}(t,\varphi_{\omega}(s,t,\eta))x,\quad
 x(s) = Id_X.
 \end{equation}
Combining \eqref{dH} and \eqref{dvarphi} we obtain
\begin{eqnarray}\label{derivetiveG}
D_2G_{\omega}(t,\eta) &=& Id_X-\int\limits_{\tau_0}^{t}\Phi_{\omega}(t,s)\left(D_1D_3\varphi_{\omega}(s,t,\eta)-A_{\omega}(t)D_3\varphi_{\omega}(s,t,\eta)\right)\,ds\nonumber\\
&=& Id_X-\int\limits_{\tau_0}^{t}D_1\left(\Phi_{\omega}(t,s)D_3\varphi_{\omega}(s,t,\eta)\right)\,ds\nonumber\\
&=& \Phi_{\omega}(t,\tau_0)D_3\varphi_{\omega}(\tau_0,t,\eta).\end{eqnarray}

From Lemma~\ref{evol}, we know that $\Phi_{\omega}(t,\tau_0)\in GL(X)$ is invertible. Moreover, \eqref{dvarphi} demonstrates that $D_3\varphi_{\omega}(\tau_0,t,\eta)$ is a transition matrix of $\dot{x} = [A_{\omega}(t)+D_2F_{\omega}(t,\varphi_{\omega}(t,s,\eta))]x$, and hence it is also invertible. Consequently, the invertibility of $D_2G_{\omega}(t,\eta)$ follows.

Hence, by  Theorem~\ref{ift}, the map $ \xi \mapsto G_{\omega}(t, \xi) $ is a local  $ C^m $ diffeomorphism for any fixed $t\in{I}$.
As the final step, since we have shown that
$G_{\omega}$
  is a homeomorphism, therefore it preserves the compactness of sets. Thus, we can apply Theorem~\ref{cacc} to
$G_{\omega},$
  which completes the proof.

\end{proof}

\section{Random dynamical systems}

In the previous section, we presented a smooth linearization result specifically for semilinear  Carath{\'e}odory differential equations in general Banach spaces. Our current objective is to apply this result to continuous time random dynamical systems, considering both random differential equations and stochastic differential equations. In order to mimic the constructions from \cite{Arnold, Wanner_95}, we restrict ourselves to finite-dimensional Banach spaces $X$. Before delving into the details, we  first introduce the research object of this section.

A random dynamical system is an object consisting of a measurable dynamical system and a cocycle over this system~\cite{Arnold}, therefore we begin with a definition of a measurable dynamical system.

\begin{definition}
Let $(\Omega, \mathcal{F},\mathbb{P})$ be a probability space. We call $(\Omega, \mathcal{F},\mathbb{P},(\theta_t)_{t\in I})$ to be a \emph{measurable dynamical system (MDS)} , if the mapping $\theta_t:\Omega\to\Omega$ satisfies the following conditions:
\begin{itemize}
    \item[(i)] the mapping $(t,\omega)\mapsto\theta_t\omega$ is measurable;
    \item[(ii)]the family $(\theta_t)_{t\in I}$ forms a group, i.e.  $\theta_0 = Id_{\Omega}$ and $\theta_{t+s} = \theta_t\circ\theta_s,$ for arbitrary $s,t\in I;$
    \item[(iii)] the mapping $\theta_t$ is $\mathbb{P}-$preserving, i.e., for arbitrary $t\in I$ and $F\in \mathcal{F},$ the identity $\mathbb{P}(\theta_t^{-1}(F)) = \mathbb{P}(F)$ holds.
\end{itemize}
\end{definition}

An MDS $(\Omega, \mathcal{F},\mathbb{P},(\theta_t)_{t\in I})$ is called \emph{ergodic} if every $\theta_t-$ invariant set has probability $0$ or $1,$ i.e. if for all $F\in \mathcal{F}$ satisfying $\theta_t^{-1}(F) = F$ for every $t\in I$ we have either $\mathbb{P}(F) = 0$ or $\mathbb{P}(F) = 1.$
\begin{definition}

A \emph{measurable random dynamical system (RDS)} on $X$
over some given MDS $(\Omega, \mathcal{F},\mathbb{P},(\theta_t)_{t\in I})$ is a measurable mapping
$\psi : I \times \Omega \times X \to X,$ $(s, \omega, x) \to \psi(s, \omega, x)$ forming a \emph{cocycle} over $\theta_t,$ i.e.,
 the mapping $\psi(t, \omega) := \psi(t, \omega, \cdot)$ satisfies

$$\psi(0, \omega) = Id_X,$$
$$\psi(s + t, \omega) = \psi(t, \theta_s \omega) \circ \psi(s, \omega), $$
for arbitrary $s,t \in I,$ and $\omega\in \Omega.$
\end{definition}

The random dynamical system $\psi$ is called \emph{linear,} if the mapping $\psi(t,\omega)$ is linear,  \emph{continuous,} if the mapping $\psi(\cdot,\omega)$ is continuous. Furthermore, it is called \emph{smooth of class $C^m$,} if it is continuous, the mapping $\psi(t,\omega)$ is of class $C^m$ and the derivatives are continuous with respect to the $(t,x),$ for arbitrary $t\in I,$ $\omega\in\Omega,$ and $x\in X.$

\begin{definition} Let $\Phi:I\times\Omega\to X$ be a linear RDS over an ergodic MDS $\theta_t$.
    We define the \emph{Lyapunov exponents} of $\Phi$  for all $x\in X, x\neq 0$  as
    $$\lambda(\omega,x) = \limsup_{t\to\infty}\frac{1}{t}\norm{\Phi(t,\omega)x}.$$
We denote the distinct values that \(\lambda(\omega, x)\) can take for \(x \neq 0\) as
\[
-\infty < \lambda_{k(\omega)}(\omega) < \dots < \lambda_1(\omega).
\] The value \(\lambda_1(\omega)\)  is called the \emph{top Lyapunov exponent.}

In the setting of the Multiplicative Ergodic Theorem~\ref{MET}, the sets $U_{\lambda}:=\{x:\lambda(\omega,x)\leq\lambda(\omega)\}$ are linear subspaces of $X,$ $U_i:=U_{\lambda_i}, $ and they form a filtration (flag of subspaces)
$$\{0\}\subset U_{k(\omega)}\subset...\subset U_1 = X,$$
where
$$\lambda(\omega,x) = \lambda_i(\omega)\iff x\in U_i\setminus U_{i+1},\quad i=1,...,k(\omega).$$
We say that $d_i(\omega):=\dim U_i-\dim U_{i+1}$ is the \emph{multiplicity} of $\lambda_i(\omega).$
%The set $$S(\theta,\Phi):=\{\lambda_i(\cdot),d_i(\cdot):i=1,...,k(\cdot)\}$$
%is called the \emph{Lyapunov spectrum} of $\Phi.$
\end{definition}

In order to apply the results from the previous section to the random dynamical system, we  first answer the following questions:

\begin{itemize}
\item[(1)] What is the connection between the semilinear  Carath{\'e}odory differential equations \eqref{deq} and the corresponding random dynamical system?
\item[(2)] How should we choose the new norm depending on both time and randomness in order to achieve our goal?
\end{itemize}

The answer to the first question is given in the following Lemma.
\begin{lemma}
    Let $(\Omega,\cF,\mathbb{P},(\theta_t)_{t\in I})$ be a MDS. If $\varphi_{\omega}$ denotes the general solution of the Carath{\'e}odory differential equation \eqref{deq},   then the following statements are equivalent:
    \begin{itemize}
        \item[(i)] The mapping $\psi_{\omega}$ defined as $\psi_{\omega}(t,x):=\varphi_{\omega}(t, 0,x)$ is a measurable RDS over $(\Omega,\cF,\mathbb{P},(\theta_t)_{t\in I}).$
        \item[(ii)] The general solution of \eqref{deq} satisfies \begin{equation}\label{lem32}
            \varphi_{\omega}(t,\tau, x) = \varphi_{\theta_{\tau} \omega}(t-\tau, 0, x)
        \end{equation} for arbitrary $\tau,t\in I$ and $x\in X.$
    \end{itemize}
\end{lemma}
\begin{proof}

$\\(i)\Rightarrow(ii)$

From the cocycle property of RDS we know that \begin{equation}\label{a}
 \varphi_{\omega}(s+t,0,x)=\varphi_{\theta_{s} \omega}(t,0,\varphi_{\omega}(s,0,x)).
\end{equation}

Moreover, it is well known that the general solution of \eqref{deq} satisfies the property
\begin{equation}\label{b}
\varphi_{\omega}(t_2,t_0,x)=\varphi_{\omega}(t_2,t_1,\varphi_{\omega}(t_1,t_0,x)).
\end{equation}

Thus,
\begin{equation*}\varphi_{\omega}(t,\tau, x) = \varphi_{\omega}((t-\tau)+\tau,\tau, x)\stackrel{\eqref{a}}{=}\varphi_{\theta_{\tau} \omega}(t-\tau,\tau,\varphi_{\omega}(\tau,0,x)) \stackrel{\eqref{b}}{=}\varphi_{\theta_{\tau} \omega}(t-\tau,0,x).\end{equation*}

$(ii)\Rightarrow(i)$

In order to prove this implication it is enough to check the cocycle property.

Namely,
\begin{equation*}\varphi_{\omega}(s+t,0, x) \stackrel{\eqref{b}}{=} \varphi_{\omega}(s+t,s,\varphi_{\omega}(s,0,x)) \stackrel{\eqref{lem32}}{=}\varphi_{\theta_{s} \omega}(t,0,\varphi_{\omega}(s,0,x)),\end{equation*}
so the lemma follows.
\end{proof}
To address the second question, we introduce a new norm in the following way.

In the setting of the Multiplicative Ergodic Theorem~\ref{MET}, suppose that $\Phi$ is a linear RDS with the Lyapunov exponents $\lambda_1(\omega),...,\lambda_{k(\omega)}(\omega).$
It is well known from, e.g., \cite{Arnold}, that when $\theta_t$ is ergodic, $k(\omega)$ and $\lambda_i(\omega)$ are independent of $\omega$. We assume ergodicity  for the remainder of this
paper.

Let $a>0$ denote a fixed real constant such that the intervals $[\lambda_i -a, \lambda_i + a],$ $i=1,...,k,$ are disjoint. For $\omega\in\Tilde{\Omega}$ and $x=x^1+...+x^{k(\omega)}\in U_1(\omega)\oplus...\oplus U_{k(\omega)}(\omega)$, we define
 a new norm \begin{equation*}|x|_{\omega} = \sqrt{|x^1|^2_{\omega}+...+|x^{k(\omega)}|^2_{\omega}} \quad,\end{equation*} where $\displaystyle|x^i|_{\omega} = \left(\int_{0}^{\infty}\|\Phi(t,\omega)x^i\|^2 e^{-2(\lambda_i t+at)}dt\right)^\frac{1}{2}$ for $x^i\in U_i(\omega).$  For $\omega\notin \Tilde{\Omega},$ we set $|x|_{\omega} = \|x\|.$  Then (according to \cite{Arnold, Wanner_95}) the following properties holds
 \begin{itemize}
     \item[(i)] $|\cdot|_{\omega}$ is a random norm on $X;$
     \item[(ii)] for every $\varepsilon >0,$ there is a measurable mapping $B_{\varepsilon}:\Omega \to [1,\infty)$ such that, for every $x\in X,$ and almost every $\omega \in \Tilde{\Omega},$ and $t \in I,$
     \begin{align*}
		\frac{1}{B_{\varepsilon}(\omega)}\|x\|&\leq|x|_{\omega}\leq B_{\varepsilon}(\omega)\|x\|,&\\
		B_{\varepsilon}(\omega)e^{-\varepsilon t}&\leq B_{\varepsilon}(\theta_t\omega)\leq B_{\varepsilon}(\omega)e^{\varepsilon t};
	\end{align*}
     \item[(iii)] for almost every $\omega \in \Tilde{\Omega},$ $t \in I,$ and every $i=1,...,k(\omega), x\in U_i(\omega),$
     \begin{equation*}
         e^{(\lambda_i-a)t} \leq \left|\Phi(t,\omega)|_{U_i(\omega)}\right|_{\omega,\theta_t\omega}\leq e^{(\lambda_i+a)t},
     \end{equation*}
     where $\left|\Phi(t,\omega)|_{U_i(\omega)}\right|_{\omega,\theta_t\omega}:=\sup\{\left|\Phi(t,\omega)x\right|_{\theta_t\omega}:x\in U_i(\omega),|x|_{\omega}\leq 1\}.$
 \end{itemize}

 It is important to note that the constant
$a$ is chosen in such a way that
$\lambda_i+a$
remains a negative expression for all
$i=1,...,k.$
In the following,  we will focus on the case, when all of Lyapunov exponents of the system are negative, i.e., $\lambda_{k}<...<\lambda_1<0$ .

%In the context of random dynamical systems, the Lyapunov exponent measures the average exponential rate of divergence or convergence of nearby trajectories. When all the Lyapunov exponents of a random dynamical system are negative, it indicates that the system expose stable behavior over time.

\begin{remark}\label{remark}
\emph{Let $\Phi$ be a linear RDS satisfying the Multiplicative Ergodic Theorem~\ref{MET} and suppose  $\Phi_{\omega}(t,s)$ is the evolution operator of a linear random differential equation $\dot{x} = A_{\omega}(t)x.$ In this case  $(H1)$ holds with $\norm{x}_{t,\omega}:=|x|_{\theta_t\omega},$ $K(\omega) = 1$ and $\alpha = -\lambda_1-a.$}

Indeed, in \cite{Arnold, Wanner_95} is shown that $\Phi$ and $\Phi_{\omega}(t,s)$ are related by $$\Phi_{\omega}(t,s) = \Phi(t-s,\theta_s\omega). $$

 From estimates in (iii) for $\left|\Phi(t,\omega)\right|_{\omega,\theta_t\omega},$ we have
    \begin{eqnarray*}
     \left|\Phi(t,\omega)x\right|_{\theta_{t}\omega}& =& \sqrt{\left|\Phi(t,\omega)x^1\right|_{\theta_{t}\omega}+...+\left|\Phi(t,\omega)x^k\right|_{\theta_{t}\omega}} \nonumber \\
     &\leq& \sqrt{(e^{(\lambda_1+a)t}|x^1|_{\omega})^2+...+(e^{(\lambda_k+a)t}|x^k|_{\omega})^2}\nonumber\\
&\leq& e^{(\lambda_1+a)t}\sqrt{|x^1|_{\omega}^2+...+|x^k|_{\omega}^2}\nonumber\\
&\leq& e^{(\lambda_1+a)t}|x|_{\omega},
    \end{eqnarray*}
and therefore for all $t\in I$
\begin{equation}\label{analog_hyp11}
    |\Phi(t,\omega)|_{\omega,\theta_t\omega} \leq e^{(\lambda_1+a)t}.
\end{equation}

   Thus,
   \begin{equation*}
   \norm{\Phi_{\omega}(t,s)}_{s,t,\omega} = |\Phi(t-s,\theta_s\omega)|_{\theta_s\omega,\theta_t\omega} = |\Phi(t-s,\theta_s\omega)|_{\theta_s\omega,\theta_{t-s}\theta_s\omega},
   \end{equation*}
   which together with~\eqref{analog_hyp11} gives us
   \begin{equation}\label{Psi_theta_S}
       \norm{\Phi_{\omega}(t,s)}_{s,t,\omega}\leq e^{(\lambda_1+a)(t-s)}.
   \end{equation}
In the above definition of the family of norms $\norm{\cdot}_{t,\omega},$ the function $\ell$ can be chosen as $\ell(t,\omega):=B_1(\theta_t\omega),$ since $B_1(\theta_{\cdot}\omega)$ is locally bounded for any $\omega\in\Tilde{\Omega}.$
\end{remark}

\subsection{$C^m-$ linearization of  random dynamical systems}

After establishing the connection between random dynamical systems and  Carath{\'e}odory differential equations, we  apply the results of the second section to RDS.  To this  end, let $\Phi$ be a linear random dynamical system on $X$ over the ergodic MDS $(\Omega,\cF, \mathbb{P}, (\theta_t)_{t\in I}),$ satisfying conditions of  the MET and generated by the following linear equation
\begin{equation}\label{rds2}
    \dot{x} = A(\theta_t\omega)x.
\end{equation}
Together with $\Phi$ we will consider
a random dynamical system $\psi$ on $X$ over the ergodic MDS $(\Omega,\cF, \mathbb{P}, (\theta_t)_{t\in I}),$ generated by
\begin{equation}\label{rds1}
    \dot{x} = A(\theta_t\omega)x+ F(\theta_t\omega,x).
\end{equation}
Assume that $\psi$ has $0$ as a fixed point, i.e. $\psi(t,\omega,0)= 0$ for almost all $t\in I,$ $\omega\in\Omega.$
Let the "nonlinear" part of $\psi$ be given by $\Psi(t,\omega,x) = \psi(t,\omega,x)-\Phi(t,\omega)x.$

While the validity of $(H1)$ has already been shown, to satisfy the remaining hypotheses, it is necessary to impose additional conditions on the random dynamical system.

\begin{enumerate}
\item[$(H_3^0)$] Assume that for almost all $t\in I,$ $\omega\in\Omega$ and all $x,\bar{x}\in X,$ there exist reals $L,M\geq 0$ such that
	\begin{align}
		|\Psi(t,\omega,x)|_{\theta_t\omega}&\leq M
		\label{hyp30}\\
		|\Psi(t,\omega,x)-\Psi(t,\omega,\bar x)|_{\theta_t\omega}&\leq L|x-\bar x|_{\omega}.
		\label{hyp31}
	\end{align}
\end{enumerate}

Then the following theorem establishes the connection between a random dynamical system $\psi$ and its linearization.
\begin{theorem}[Topological linearization of RDS]\label{thm_top}
    Assume that the linear random dynamical system $\Phi,$ generated by \eqref{rds2},  satisfies the assumptions of the MET \ref{MET} and has negative Lyapunov exponents. Assume that
   the nonlinear part $\Psi$  of a RDS $\psi,$ generated by \eqref{rds1},  satisfies $(H_3^0).$ If $L\leq \alpha,$
    then the RDS $\psi $ and $\Phi$ are topologically equivalent, i.e. there is a measurable mapping $h:\Omega \times X\to X$ with the following properties:
    \begin{enumerate}
        \item[(i)]  $h(\omega) = h(\omega,\cdot)$ is a homeomorphism on $X$ with $h(\omega, 0) = 0$  for almost all $\omega\in\Omega;$
        \item[(ii)] $h(\omega)$ maps $\omega-$orbits of $\Phi$ onto $\omega-$orbits of $\psi$ in the sense that for every $\omega\in\Omega,$ arbitrary  $t\in I$ and $x\in X$
        $$h(\omega)\circ\Phi(t,\omega)x = \psi(t,\omega,h(\omega, x));$$
        \item[(iii)] h has property of being near identity, i.e. for all $\xi,\eta\in X$
        $$|h(\omega,\xi)-\xi|_{\theta_t\omega}\leq\frac{M}{\alpha} \quad\text{and}\quad |h^{-1}(\omega, \eta)-\eta|_{\theta_t\omega}\leq\frac{M}{\alpha},$$
        where $\alpha = -\lambda_1-a > 0.$
    \end{enumerate}
\end{theorem}

\begin{proof}
Applying Theorem \ref{thm_lin} $\omega$-by-$\omega$ to the system \eqref{rds1} with $\Omega$ replaced by $\Tilde{\Omega}$ will yield a mapping $h:\Omega\times I\times X\to X,$ defined as follows
\begin{equation}\label{h(w)}
    \begin{cases}h(\omega,\xi) = \xi+\int_{0}^{t}\Phi(t-s,\theta_s\omega)\Psi(s,\theta_s\omega,\varphi_{\theta_t\omega}^*(\xi)(s)+\Phi(t-s,\theta_s\omega)\xi)\d s,\quad \omega\in\Tilde{\Omega},\\
   h(\omega,\xi) = 0,\quad\omega\notin\Tilde{\Omega},
   \end{cases}
\end{equation}
with the inverse
\begin{equation}\label{h-1(w)}
    \begin{cases}h^{-1}(\omega,\eta) = \eta-\int_{0}^{t}\Phi(t-s,\theta_s\omega)\Psi(s,\theta_s\omega,\psi(s,\omega,\eta))\d s,\quad \omega\in\Tilde{\Omega},\\
   h(\omega,\eta) = 0,\quad\omega\notin\Tilde{\Omega},    \end{cases}
\end{equation}
such that statements $(i)$ and $(ii)$ holds.

Let us check that the homeomorphism defined in \eqref{h(w)} and its inverse are near identity. Indeed,
it is easy to verify that \begin{eqnarray*}
   |h^{-1}(\omega,\eta)- \eta|_{\theta_t\omega}&\stackrel{\eqref{h-1(w)}}{\leq}&\int_{0}^{t}|\Phi(t-s,\theta_{s}\omega)|_{\theta_s\omega,\theta_t\omega}|\Psi(s,\theta_s\omega,\psi(s,\omega,x))|_{\theta_t\omega}\,ds\nonumber\\
   &\stackrel{\eqref{Psi_theta_S}}{\leq}&\int_{0}^{t}e^{(\lambda_1+a)(t-s)}|\Psi(s,\theta_s\omega,\psi(s,\omega,x))|_{\theta_t\omega}\,ds\stackrel{\eqref{bound_psi}}{\leq} M \int_{0}^{t}e^{(\lambda_1+a)(t-s)}\,ds\nonumber\\
   &\leq&\frac{M}{-\lambda_1-a}.
\end{eqnarray*}
Using the same arguments, we obtain
\begin{equation*}
    |h(\omega,\xi)-\xi |_{\theta_t\omega}\leq\frac{M}{-\lambda_1-a}.
\end{equation*}

%Moreover, the fact that homeomorphism defined in \eqref{h(w)} maps $\omega-$ orbits of $\psi$ onto $\omega-$ orbits of $\Phi$ is a direct consequence of the following thoughts:
%\begin{eqnarray*}
    %    \frac{d }{d t}h(\omega)\Phi(t,\omega)\xi
   %     &=&\frac{d }{d t}\Phi(t,\omega)\xi+\frac{d }{d t}\varphi^*_{\theta_t\omega}(\xi)(t)\\
%&=&A(\theta_t\omega)(\Phi(t,\omega)\xi +  \varphi^*_{\theta_t\omega}(\xi)(t)) + F(\theta_t\omega ,h(\omega)\Phi(\omega,t)\xi) \\
   %     &=& A(\theta_t\omega)h(\omega)\Phi(t,\omega)\xi+F(\theta_t\omega, h(\omega)\Phi(t,\omega)\xi)\\
  %      &=&\psi(t,\omega,h(\omega)\xi).
  %  \end{eqnarray*}
\end{proof}
\begin{theorem}[Smooth linearization of RDS]\label{thm_smooth}
    Assume that the linear random dynamical system $\Phi$ generated by \eqref{rds2}  satisfies the assumptions of the MET \ref{MET}, has negative Lyapunov exponents. Assume further that
      all partial derivatives $D_2^{j}\Psi:I\tm\Omega\tm X\to L_j(X)$, $1\leq j\leq m$ satisfy Carath{\'e}odory property, and, moreover, there exist reals $M_j\geq 0$ such that
	\begin{align}\label{bound_psi}
		|D_2^{j}\Psi(\omega,x)|_{\theta_t\omega}&\leq M_j\quad\text{for almost all } t\in I, \omega\in\Omega,  \text{ all } x\in X,  0\leq j\leq m,
	\end{align}
 and  $M_1\leq \alpha,$  $\alpha = -\lambda_1-a, $ $ \alpha>0, $
 then all statements of the previous theorem hold and $h(\omega,\cdot)$ is a $C^m-$ diffeomorphism.
\end{theorem}
\begin{proof}
    It suffices to demonstrate that the homeomorphism~\eqref{h(w)} is $C^m$. To do this we express the $D_2h(\omega,\cdot)^{-1}$ in a direct form, using the same arguments as in Theorem \ref{thmdiff} :
    $$D_2h(\omega,\cdot)^{-1}= \Phi(t,\omega)D_3\psi(s,\omega,\cdot).$$

  Applying Theorem \cite[Thm. 2.2.2, p.~60]{Arnold},  we see that
   the determinant $\det D\psi(t,\omega,x)$ satisfies Liouville's equation:
        $$\det D\psi(t,\omega,x) = \exp \int_{0}^{t}\text{ trace } D\Psi(\theta_s\omega,\psi(s,\omega,x))\d s.$$
  Moreover, the linear part of RDS satisfies:
    $$\det \Phi(t,\omega) = \exp\int_{0}^{t} \text{ trace } \Phi(t-s,\theta_s\omega)\d s.$$ Therefore, $\det(\Phi(t,\omega)D_3\psi(s,\omega,\eta))>0$ for any $t\in I,$ $\omega\in\Omega$ which gives us the invertibility of $D_2h^{-1}(\omega,\eta)$ and the theorem follows.
\end{proof}

The assurance of the existence and smallness of a Lipschitz constant, as required in the last theorem, is, of course, not guaranteed in general. Therefore, we will adjust the nonlinear part of the vector fields by cutting them off when they exceed a certain threshold. This will bring us to the local linearization result, which is more suitable to stochastic systems and can be applied in the next section.

For deriving local results from global ones we  consider a random dynamical system  $\psi$ having a fixed point at the origin. Our objective is to construct a novel random dynamical system, denoted as $\Tilde{\psi}$, which faithfully reproduces the behavior of $\psi$ within a localized random neighborhood of origin.  To make this construction easier, we use the insights from the following lemmas.

\begin{lemma}\label{wan_1}\cite[Lem. 4.2, p. 262]{Wanner_95}
   Let $L_0 > 0,$ $(\Omega,\cF, \mathbb{P}, (\theta_t)_{t\in I})$ denote a measurable dynamical system, and consider a measurable mapping $F:\Omega\times X\to X$ such that $F(\omega,\cdot)$ is continuous with $F(\omega,0) = 0$ for all $\omega\in\Omega.$ If
   \begin{equation}\label{lem_cut_of1_1}
     \lim_{(x,y)\to(0,0)}\frac{\norm{F(\omega,x)-F(\omega,y)}}{\norm{x-y}} = 0
   \end{equation}
   and
   \begin{equation}\label{lem_cut_of1_2}
     \sup_{\|x\|\leq c, \|y\|\leq c, x\neq y}\frac{\norm{F(\omega,x)-F(\omega,y)}}{\norm{x-y}} \in L^1(\Omega,\cF, \mathbb{P})
   \end{equation}
   hold for some $c>0,$  then there are measurable mappings $\sigma:\Omega\to (0,c],$ $L:\Omega\to \mathbb{R}^+$ and $\Tilde{F}:\Omega\times X\to X,$ as well as a $\theta_t-$ invariant set $\hat{\Omega}\in \cF$ with $\mathbb{P}(\hat{\Omega}) = 1$ such that the following holds:
   \begin{itemize}
    \item[(i)] if we define the random neighborhood of $0$ by
    \begin{equation}\label{U}
    U(\omega):=\{x\in X:\|x\|< \sigma(\omega)\},\end{equation}
    then the identity $F(\omega,x) = \Tilde{F}(\omega,x)$ holds for all $\omega\in\hat{\Omega}$ and $x\in U(\omega).$ Additionally, we have $\Tilde{F}(\omega,x) = 0$ for $\omega\in\Omega\setminus\hat{\Omega}$ and $x\in X;$
    \item[(ii)] for arbitrary $\omega\in\hat{\Omega}$ and $x,y\in X$ we have
    $$|\Tilde{F}(\omega,x)-\Tilde{F}(\omega,y)|_{\omega}\leq L(\omega)|x-y|_{\omega}\quad\text{ and }\quad |\Tilde{F}(\omega,x)|_{\omega}\leq L(\omega)c,$$
    as well as $\displaystyle\int_{0}^{1}L(\theta_s\omega)ds\leq L_0;$
    \item[(iii)] for all $\omega\in\hat{\Omega}$ and $a,b\in \mathbb{R}$ with $a\leq b$ the estimate $\displaystyle\inf_{a\leq t\leq b}\sigma(\theta_t\omega)>0$ holds, i.e., the mapping $\sigma(\theta_t\omega)$ is locally bounded away from $0.$
   \end{itemize}
\end{lemma}

\begin{lemma}\cite[Prop. 4.3, p. 264]{Wanner_95}\label{lema_wanner}
    Assume we are given a continuous random dynamical system $\psi$ on $X$ over the ergodic measurable dynamical system $(\Omega,\cF, \mathbb{P}, (\theta_t)_{t\in I})$ with fixed point $0,$ which is generated by the random differential equation \eqref{rds1}.
If $\norm{A}\in L^1(\Omega,\cF, \mathbb{P}),$ $F(\omega,0) = 0$ and $$\lim_{(x,y)\to(0,0)}\frac{\norm{F(\omega,x)-F(\omega,y)}}{\norm{x-y}} = 0 $$
for all $\omega\in \Omega,$
then the following holds
\begin{itemize}
    \item[(i)] There is a $\theta_t-$ invariant set $\hat{\Omega}\in\cF$ with $\mathbb{P}(\hat{\Omega}) = 1,$ as well as a measurable mapping $\sigma:\Omega\to\mathbb{R}^+$ and a continuous random dynamical system $\Tilde{\psi}(t,\omega,x) = \Phi(t,\omega)x+\Tilde{\Psi}(t,\omega,x)$ over $(\Omega,\cF, \mathbb{P}, (\theta_t)_{t\in I})$ satisfying $(H_3^0),$ as well as
\begin{equation}\label{loc_lin_eq}
  \psi(t,\omega,x) = \Tilde{\psi}(t,\omega,x)
\end{equation}
for $\omega\in\hat{\Omega},$ $x\in U(\omega)$ and $t\in I_{\max}(\omega,x):=\{0,t_{\max}(\omega,x)\},$ where $U(\omega)$ is defined in \eqref{U} and
\begin{equation}\label{t_max}
t_{\max}(\omega,x):=\sup\{t\in I: \psi(\tau,\omega,x)\in U(\theta_{\tau}\omega) \fall 0\leq\tau\leq t\}\geq 0. \end{equation}
\item[(ii)] For every $\omega\in\hat{\Omega}$ we have $\displaystyle\lim_{x\to 0} t_{\max}(\omega,x) = \infty.$ In particular, this implies that for every choice of $\omega\in\hat{\Omega}$ and $\tau_0\in I$ the existence of a neighborhood $U_{\tau_0}(\omega)\subset X$ of the origin such that the equality \eqref{loc_lin_eq} holds for all $t\in (0,\tau_0)$ and $x\in U_{\tau_0}(\omega).$
\end{itemize}

\end{lemma}

\begin{lemma}
   If the conditions of  Lemma~\eqref{lema_wanner} holds, then the nonlinear part $\Tilde{\Psi}$ of the  newly defined  local random dynamical system $\Tilde{\psi}$ satisfies the assumptions contained in $(H^0_3)$ with $\Omega$ replaced by $\hat{\Omega}.$
\end{lemma}
\begin{proof}
Let $\tilde{\Phi}_{\omega}(t,s)$ denote the evolution operator of the linear random differential equation $\dot{x} = A(\theta_t\omega)x.$ The correspondence between the linear random dynamical system $\Phi(t-s,\theta_s\omega)$ and $\tilde{\Phi}_{\omega}(t,s)$  shown in Remark \ref{remark}, together with \eqref{analog_hyp11} implies \begin{equation}\label{loc_1}
|\tilde{\Phi}_{\omega}(t,s)|_{\theta_s\omega,\theta_t\omega} = |\Phi(t-s,\omega)|_{\theta_s\omega,\theta_t\omega} \leq e^{(\lambda_1+a)(t-s)}
\end{equation}
$\text{ for arbitrary } s\leq t, \omega\in\hat{\Omega}.$

Let $\Tilde{F}$ denote the mapping guaranteed by Lemma \ref{wan_1}. Since $\Tilde{\psi}(t,\omega,x)$ solves the initial value problem $\dot{x} = A(\theta_t\omega)x+\Tilde{F}(\theta_t\omega,x),$ $x(0) = x,$ from the variation of constants formula we have
$$\Tilde{\psi}(t,\omega,x) = \Tilde{\Phi}_{\omega}(t,0)x+\int_0^t\Tilde{\Phi}_{\omega}(t,s)\Tilde{F}(\theta_s\omega,\Tilde{\psi}(s,\omega,x))ds,$$
and therefore
$$\Tilde{\Psi}(t,\omega,x) = \int_0^t\Tilde{\Phi}_{\omega}(t,s)\Tilde{F}(\theta_s\omega,\Tilde{\Phi}_{\omega}(s,0)x+\Tilde{\Psi}(s,\omega,x))ds,$$
for arbitrary $t\in I,$ $x\in X$ and $\omega\in\hat{\Omega}.$
This implies the estimate
\begin{eqnarray*}
    |\Tilde{\Psi}(t,\omega,x)|_{\theta_t\omega}&\leq&\int_0^t|\Tilde{\Phi}_\omega(t,s)|_{\theta_s\omega,\theta_t\omega}|\Tilde{F}(\theta_s\omega,\Tilde{\Phi}_{\omega}(s,0)x+\Tilde{\Psi}(s,\omega,x))|_{\theta_s\omega}\,ds
    \\&\leq& c \int_0^t e^{(\lambda_1+a)(t-s)}L(\theta_s\omega)\,ds\leq c e^{(\lambda_1+a)}\int_0^1L(\theta_s\omega)ds\\&\leq& L_0 c e^{(\lambda_1+a)},
\end{eqnarray*}
for any $t\in [0,1],$ $\omega\in\hat{\Omega}$ and $x\in X.$

Now, our aim is to verify the Lipschitz condition. For any $t\in [0,1],$ $\omega\in\hat{\Omega}$ and $x,\bar{x}\in X$ we have:
\begin{eqnarray*}
    |\Tilde{\Psi}(t,\omega,x)-\Tilde{\Psi}(t,\omega,\bar{x})|_{\theta_t\omega}&\leq&\int_0^t|\Tilde{\Phi}_\omega(t,s)|_{\theta_s\omega,\theta_t\omega}|\Tilde{F}(\theta_s\omega,\Tilde{\Phi}_{\omega}(s,0)x+\Tilde{\Psi}(s,\omega,x))\\&-&\Tilde{F}(\theta_s\omega,\Tilde{\Phi}_{\omega}(s,0)\bar{x}+\Tilde{\Psi}(s,\omega,\bar{x}))|_{\theta_s\omega}\,ds
    \\&\leq& c \int_0^t e^{(\lambda_1+a)(t-s)}L(\theta_s\omega)(|\Tilde{\Phi}_{\omega}(s,0)(x-\bar{x})|_{\theta_s\omega}\,ds\\
    &+&|\Tilde{\Psi}(s,\omega,x)-\Tilde{\Psi}(s,\omega,\bar{x})|_{\theta_s\omega})\,ds.
\end{eqnarray*}
Denote by $\vartheta(t) := e^{-(\lambda_1+a)t}|\Tilde{\Psi}(s,\omega,x)-\Tilde{\Psi}(s,\omega,\bar{x})|_{\theta_s\omega}.$ Then the last estimation can be expressed as follows
\begin{eqnarray*}
  \vartheta(t)&\leq&\int_0^t e^{-(\lambda_1+a)s} L(\theta_s\omega)|\Tilde{\Phi}_{\omega}(s,0)|_{\omega,\theta_s\omega}|x-\bar{x}|_{\omega}\,ds+ \int_0^t L(\theta_s\omega)\vartheta(s)\,ds\\
  &\leq&|x-\bar{x}|_{\omega}\int_0^tL(\theta_s\omega)\,ds + \int_0^t L(\theta_s\omega)\vartheta(s)\,ds.
\end{eqnarray*}
By applying Gronwall's inequality, we obtain
\begin{eqnarray*}
  \vartheta(t)&\leq& |x-\bar{x}|_{\omega}\int_0^tL(\theta_s\omega)\,ds + \int_0^t |x-\bar{x}|_{\omega}L(\theta_s\omega)\left(\int_0^s L(\theta_\tau\omega)\,d\tau\right)e^{\int_s^t L(\theta_\tau\omega)\,d\tau}\,ds\\
  &\leq& |x-\bar{x}|_{\omega}(L_0+L_0^2e^{L_0}),
\end{eqnarray*}
for any $s,t\in[0,1],s\leq t,$ $x,\bar{x}\in X$ and $\omega\in\hat{\Omega}.$

 Finally, we obtain
 $$|\Tilde{\Psi}(t,\omega,x)-\Tilde{\Psi}(t,\omega,\bar{x})|_{\theta_t\omega}\leq e^{\lambda_1+a}(L_0+L_0^2e^{L_0})|x-\bar{x}|_{\omega},$$
 which completes the proof.
 \end{proof}

 Now, we are able to formulate a local $C^m$  result  as well.

\begin{theorem}
  [Local topological linearization of RDS]\label{thm_local}
Let $\psi$ be a  given  continuous random dynamical system on $X$ over the ergodic metric dynamical system $(\Omega,\cF, \mathbb{P}, (\theta_t)_{t\in I})$ with fixed point $0,$  generated by the random differential equation \eqref{rds1}. And and $\Phi$ be a linear random dynamical system  on $X,$ generated by the random differential equation \eqref{rds2}.
If $\Phi$ satisfies MET, has negative Lyapunov exponents, $F(\omega,0) = 0$ and $$\lim_{(x,y)\to(0,0)}\frac{\norm{F(\omega,x)-F(\omega,y)}}{\norm{x-y}} = 0 $$
for all $\omega\in \Omega,$
then the  newly defined  local random dynamical system $\Tilde{\psi}$ and $\Phi$ are topologically equivalent on $I_{\max}(\omega,x):=(0,t_{\max})$ in the sense of Theorem \ref{thm_top} for all $x\in U(\omega)$ and $\omega\in\hat{\Omega},$ where $t_{\max}$ defined in \eqref{t_max} and $U(\omega)$ defined in \eqref{U}.
\end{theorem}

\begin{theorem}[Local smooth linearization of RDS]\label{thm_local_smooth}
    If the assumptions of Theorem~\eqref{thm_local} hold and additionally there exist reals $M_j\geq 0$ such that all partial derivatives $D_2^{j}F:\hat{\Omega}\tm U(\omega)\to L_j(U(\omega))$, $1\leq j\leq m$ satisfy
	\begin{align}
		|D_2^{j}F(\omega,x)|_{\theta_t\omega}&\leq M_j\quad\text{for almost all $t\in I_{\max}$, $\omega\in\hat{\Omega},$  all $x\in U(\omega)$ and }0\leq j\leq m,
	\end{align}
 then random dynamical system $\Tilde{\psi}$ and $\Phi$ are $C^m-$topologically equivalent  on $I_{\max}.$
\end{theorem}
\begin{remark}\emph{
     The proof of both above theorems \ref{thm_local} and \ref{thm_local_smooth} follow from the application of the Theorems \ref{thm_top} and \ref{thm_smooth} respectively to the new random dynamical system $\Tilde{\psi}$ with $X$ restricted to $U(\omega)$ and $I$ restricted to $I_{\max}(\omega,x).$ }
\end{remark}

\subsection{$C^m-$ linearization of stochastic differential equations}

In this section, we establish a local correspondence between a random dynamical system, generated by a stochastic differential equation, and its linearization by using conjugate transformations and applying previously established local linearization results.  Assume that $X$ is a finite dimensional Banach space. We will also need the following notations borrowed from \cite{Kunita}. For $k\in \mathbb{N}, \delta\in(0,1],$ denote by $C^{k,\delta}_{b}$ the set of functions $f:X\to X$ which have partial derivatives $D^{i}f$ of order up to $|\alpha|\leq k,$ are linearly bounded and for which the derivatives of order $k$ are $\delta-$H\"older continuous. In other words,
\begin{equation*}
    \sup_{x\in X}\frac{|f(x)|}{1+|x|}+\sum\limits_{|\alpha|=1}^{k}\sup_{x\in X}|D^{\alpha}f(x)|+\sum_{|\alpha| = k}\sup_{\substack{x, \bar{x} \in X \\ x \neq \bar{x}}}\frac{|D^{\alpha}f(x)-D^{\alpha}f(\bar{x})|}{|x-\bar{x}|^{\delta}}<\infty.
\end{equation*}

In order to obtain a SDE result, the approach of \cite{imkeller} plays a key role for us. Specifically, the next two theorems demonstrate the feasibility of a suitable random stationary coordinate transformation (cohomology) that establishes a connection between cocycles of stochastic and random differential equations.
\begin{theorem}\cite[Thm. 3.1, p. 140]{imkeller}\label{il_1}
    If  $f_1,...,f_k\in C^{2,\delta}_{b}$ satisfy $\sum\limits_{i=1}^k f_{i}^j\frac{\partial f_i}{\partial x_j}\in C^{2,\delta}_{b},$ then there exists a random flow of diffeomorphisms $\Phi:I\times X\times\Omega\to X$ such that the following properties hold:
    \begin{itemize}
        \item [(i)] for any $(\eta,\tau)\in X\times I$ the process $\Phi(x,\tau)$  satisfies the stochastic integral equation
    \begin{equation}\label{imk}
        \Phi_t(x,\tau) = x+e^{-\tau}\sum_{i=1}^k\int\limits_{-\infty}^{t}e^s f_i(\Phi_s(x,\tau))\circ dW_s^i;
    \end{equation}
    \item[(ii)] if
    \begin{equation}
        H_t = \Phi_t(\cdot,t),
    \end{equation}
    \begin{equation}
        \Gamma_t = \frac{\partial}{\partial \tau}\Phi_t(\cdot,t),
\end{equation}
then $H$ is a stationary cocycle of diffeomorphisms, $\Gamma$ is a stationary random field, and for any $x\in X$  the processes $H(x)$ and $\Gamma(x)$ satisfy the SDE
\begin{equation*}
    dH_t(x) = \sum_{i=1}^k f_i(H_t(x))\circ dW_t^i + \Gamma_t(x) dt;
\end{equation*}
\item[(iii)] for $x\in X,$ $\Gamma(x)$ satisfies the stochastic integral equation
\begin{equation*}
    \Gamma_t(x) = -(H_t(x)-x)+e^{-t}\sum_{i=1}^k\int\limits_{-\infty}^t e^s\frac{\partial }{\partial x}f_i(H_t(x))\Gamma_s(x)\circ dW_s^i.
\end{equation*}
\end{itemize}
\end{theorem}
\begin{theorem}\cite[Thm. 3.2, p. 142]{imkeller}\label{il_2}
 If $f_0\in C^{1,\delta}_{b}$ and $f_1,...f_k\in C^{2,\delta}_{b}$ satisfy $\sum\limits_{i=1}^k f_{i}^j\frac{\partial f_i}{\partial x_j}\in C^{2,\delta}_{b},$ then the following holds for the random vector field $g:\Omega\times X\to X,$
 \begin{equation*}
     g(\cdot,y) := \frac{\partial}{\partial y} H_0^{-1}(y)[f_0(H_0(y))+\Gamma_0(y)],
 \end{equation*}
  with $H$ and $\Gamma$ given by the previous theorem:
 \begin{itemize}
     \item[(i)]the random differential equation \begin{equation*}
         dy_t = g(\theta_t\cdot,y_t)dt
     \end{equation*}
     generates a random cocycle of diffeomorphism $\Psi$;
      \item[(ii)] the SDE
      \begin{equation*}
          dx_t = f_0(x_t)dt+\sum_{i=1}^k f_i(x_t)\circ dW^i_t
      \end{equation*}
      generates a random cocycle of diffeomorphism $\Phi$;
   \item[(iii)] $\Phi$ and $\Psi$ are conjugate with cohomology $H_0,$ i.e., for $t\in I,$ $\omega\in\Omega$
 \begin{equation*}
     \Phi_t(\omega) = H_0(\theta_t\omega,\cdot)\circ \Psi(\omega)\circ H_0^{-1}(\omega,\cdot).
 \end{equation*}\end{itemize}

\end{theorem}

\begin{remark}
    \emph{It is important to note the following. Consider  Theorems~\ref{il_1} and \ref{il_2} in the context of a higher degree of smoothness for the vector fields  $f_1,...,f_k.$ Specifically,  suppose $f_1,...,f_k\in C^{m+1,\delta}_{b},$ and $\sum\limits_{i=1}^m f_{i}^j\frac{\partial f_i}{\partial x_j}\in C^{m+1,\delta}_{b}.$ Then $\Phi$ and, consequently, $H$ are $C^m$ diffeomorphisms.}

This fact is a direct consequence of applying the following Theorem \ref{cm_arn} to system \eqref{imk}.
\end{remark}
\begin{theorem}\label{cm_arn}\cite[Thm. 2.3.32, p. 93]{Arnold}
If $f_0\in C^{m,\delta}_{b}$ and $f_1,...,f_k\in C^{m+1,\delta}_{b}$ satisfy $\sum\limits_{i=1}^k f_{i}^j\frac{\partial f_i}{\partial x_j}\in C^{m,\delta}_{b},$
then the classical Stratonovich stochastic differential equation
\begin{equation*}
   dx_t = f_0(x_t)dt+\sum\limits_{j=1}^k f_j(x_t)\circ dW_t^j
\end{equation*}
generates a unique $C^m$ random dynamical system over the filtered dynamical system describing a Brownian motion.

\end{theorem}

Transitioning between SDEs and RDEs through  coordinate changes must maintain the principles of a multiplicative ergodic theory, such as Lyapunov exponents and Oseledets spaces. The following proposition demonstrates that these properties are indeed preserved, given that the coordinate change fulfills specific integrability criteria.

\begin{lemma}\label{lemma_imk}\cite[Prop. 5.1, p. 156]{imkeller}
Let $\Phi:I\times \Omega \to X$ be a linear cocycle such that \begin{align*}
		\sup_{0\leq t\leq 1}\log^+\norm{\Phi(t,\omega)}&\in L^1(\Omega,\cF,\mathbb{P}),&
		\sup_{0\leq t\leq 1}\log^+\norm{\Phi^{-1}(t,\omega)}&\in L^1(\Omega,\cF,\mathbb{P}).
	\end{align*}
 If $H:\Omega\times X \to X$ is a random linear mapping satisfying
 \begin{itemize}
     \item[(i)] $t\mapsto H(\theta_t\omega)$ is continuous for every $\omega\in\Omega,$
     \item[(ii)] one has the integrability condition
     \begin{align*}
		\sup_{0\leq t\leq 1}\log^+\norm{H(\theta_t\omega)}&\in L^1(\Omega,\cF,\mathbb{P}),&
		\sup_{0\leq t\leq 1}\log^+\norm{H^{-1}(\theta_t\omega)}&\in L^1(\Omega,\cF,\mathbb{P}),
	\end{align*}
    \end{itemize}
 then $$\Psi(t,\omega) = H(\theta_t\omega,\cdot)\circ\Psi(t,\omega)\circ H^{-1}(\omega,\cdot)$$
 for any $t\in I$ and almost all $\omega\in\Omega$ defines a linear cocycle which satisfies MET. $\Psi$ possesses the same Lyapunov exponents as $\Phi.$ If $U_1,...U_k$ are the Oseledets spaces of $\Phi,$ then the Oseletets spaces of $\Psi $ are given by $HU_1,...,HU_k.$

\end{lemma}

Following the preliminary discussions, we establish a local $C^m$-linearization result for random dynamical systems (RDS) generated by stochastic differential equations (SDEs).

\begin{theorem}[Local smooth linearization of SDEs]
   Let the vector fields $f_0\in C^{m,\delta}_{b}$  and $f_1,...,f_k\in C^{m+1,\delta}_{b}$ satisfy $f_i(0) = 0,$ $i=0,...,k,$ and additionally assume $$\displaystyle\sum\limits_{i=1}^k\sum\limits_{j=1}^d f_{i}^j\frac{\partial f_i}{\partial x_j}\in C^{m+1,\delta}_{b}.$$
 If the random dynamical system
 $\Phi$  generated by the linear stochastic differential equation
\begin{equation}\label{sde_1}
    dx_t = \frac{\partial }{\partial x}f_0(0)x_t \,dt+\sum\limits_{i=1}^k \frac{\partial }{\partial x}f_i(0)x_t\circ \,dW^i_t
\end{equation}
satisfies MET with negative Lyapunov exponents,
then it is $C^{m-1}-$ topologically equivalent to the
RDS $\psi$  generated by the  stochastic differential equation
\begin{equation}\label{sde_2}
    dx_t = f_0(x_t)\,dt+\sum\limits_{i=1}^k f_i(x_t)\circ \,dW^i_t
\end{equation}
for all $x\in X$, $t\in I$ and almost all $\omega\in \Omega.$
\end{theorem}
\begin{proof}

To prove the theorem, first, we need to conjugate the nonlinear SDE \eqref{sde_2} with the corresponding random differential equation.

In order to do this, let us construct $\Gamma$ and $H$ and for $y\in X$ define the random  field $g:\Omega\tm X\to X$  according to Theorems \ref{il_1}, \ref{il_2} by
$$g(\cdot,x) = \frac{\partial}{\partial x} H_0^{-1}(x)(f_0(H_0(y))+\Gamma_0(y)).$$

The random differential equation induced by $g$ has $0$ as a fixed point by our hypothesis and construction and is differentiable. Let $$A = \frac{\partial}{\partial x} g.$$
Denote by $\psi^0$ the RDS generated by $g,$ and by $\Phi^0$ the RDS generated by
$$dx_t = A(\theta_t\cdot)x_t \,dt.$$

By the construction, the RDS $\psi$ is conjugated to $\psi^0$ via the cohomology $H$, and the RDS $\Phi$ generated by \eqref{sde_1} is conjugated to $\Phi^0$ via $\displaystyle\frac{\partial H}{\partial x}$.

Now we turn to the application of the local linearization result (Theorem~\ref{thm_local}) to the newly obtained systems $\psi^0$ and $\Phi^0$. In order to do this, we have to check the assumptions of Theorem~\ref{thm_local}.

By the construction of $\Gamma$ and $H$, and the smoothness properties of $f_i$ for $0 \leq i \leq k$, it follows that all required properties for $g$ are satisfied. Moreover, from Lemma \ref{lemma_imk} and the construction of $H$, it follows that $\Phi^0$ satisfies the integrability condition of the MET. Therefore, by applying the local linearization theorem, we conclude that $\Phi$ and $\psi^0$ are  topologically equivalent via $C^m$ diffeomorphism $h$.

With all the groundwork laid, we now proceed to conclude that random dynamical systems $\Phi$ and $\psi$ are  topological equivalent via the $C^{m-1}-$ diffeomorphism    $$\tilde{h} = \displaystyle\frac{\partial }{\partial x}H_0(\theta_t\omega,\cdot)\circ h \circ H_0^{-1}(\omega,\cdot)$$
for all $x\in X$, $t\in I,$ and almost all $\omega\in \Omega.$

\end{proof}

\setcounter{section}{0}
\renewcommand{\thesection}{\Alph{section}}
\section{Appendix}
This appendix contains all the essential results and lemmas used throughout the paper to support the proofs of our results. It is organized into two logical parts: deterministic and random.

In the deterministic part, we start with Gronwall’s Lemma, which is fundamental for estimating the growth of solutions of differential equations.

Throughout, $X, Y$ denote Banach spaces.
\begin{lemma}(Gronwall's inequality)\cite[Lem. 2.5, p. 53]{Wanner_96}, \cite{bres_07}
Suppose we are given a nonempty interval I, a point $\tau\in I$ and let $\alpha,u:I\to\mathbb{R}$  be continuous  with $\alpha$ non-decreasing. If for fixed $\beta >0 $ the  inequality
$$u(t)\leq \alpha(t) + \int_{\tau}^{t}\beta u(s)ds$$
holds,
then
\begin{equation}\label{Gron1}
    u(t)\leq \alpha(t)e^{\beta (t-\tau)}\quad \fall t\in I.
\end{equation}
\end{lemma}
\begin{theorem}\label{var_const}(Variation of Constants Formula)\cite[Thm. 2.10, p. 58]{Wanner_96}, \cite{codd_56} Consider the system:
\begin{equation}\label{1}
   \dot{x} = A_{\omega}(t)x + F_{\omega}(t,x),
\end{equation}
where \( A_{\omega}: I \to L(X) \) and \( F_{\omega}: I \times X \to X \) are locally integrable.
Then the solution \( \varphi(t,s,\xi) \) of \eqref{1} is given by the variation of constants formula:
\[
\varphi(t,s,\xi) = \Phi_{\omega}(t,s)\xi + \int_{s}^t \Phi_{\omega}(t,\tau) F_{\omega}(\tau, \varphi(\tau,s,\xi)) \, d\tau,
\]
which is valid for all \( (t,s,\xi) \in I \times I \times X \).
\end{theorem}
Now we present the well-known fixed point theorem of Banach and some related results.

\begin{theorem}\label{Fixed point}(Banach's Fixed Point Theorem)\cite[Thm. B.1, p. 114]{Wanner_96}, \cite{zeidler}

If $\cT : X \to X$ a contraction, i.e.~there exists $0\leq c< 1$ such that
$$\norm{\cT(x)-\cT(\bar{x})}\leq c\norm{x-\bar{x}}\quad \text{ for all } x,\bar{x}\in X,$$
then there exists  exactly one fixed point $\varphi \in X$ of $\cT.$
\end{theorem}

\begin{theorem}\label{uniform contraction}(Uniform Contraction Principle) \cite[Thm. 2.2, p. 25]{Hale}

Suppose $U\subseteq X,$ $V\subseteq Y$ are open, let $\bar{U}$ denote the closure of $U$ and $m\in\mathbb{N}_0.$ If
  $\cT\in C^m( \bar{U}\times V,X),$  is a uniform contraction, i.e.~there exists a $0\leq c< 1$ such that
$$\norm{\cT(x,y)-\cT(\bar{x},y)}\leq c\norm{x-\bar{x}}\quad\text{ for all }x,\bar{x}\in \bar{U} \text{ and } y\in V ,$$
then the unique fixed point $\varphi(y)$ of $\cT(\cdot,y)$ in $\bar{U}$ satisfies $\varphi \in C^m(V,X).$

\end{theorem}
\begin{theorem}\label{ift}(Local Inverse Mapping Theorem) \cite[Thm. 4.F, p. 172]{zeidler}

Let $U\subseteq X, x_0\in U$ and  $f:U\to Y$ be a $C^m$ mapping for $m\in\mathbb{N}$. Then $f$ is a local $C^m-$ diffeomorphism at $x_0$ if and only if $Df(x_0)$ is bijective.
\end{theorem}

\begin{theorem}\label{cacc}(Global Inverse Mapping Theorem)\label{hadamard} \cite[Thm. 4.G, p. 174]{zeidler}

Let $f:X\to Y$ be a  local $C^m-$ diffeomorphism for $m\in\mathbb{N}$ at every point of $X.$  Then $f$ is a $C^m-$ diffeomorphism if and only if $f$ is proper.
\end{theorem}

\begin{theorem}\label{diff_int}(Differentiation under the Integral) \cite[Thm. A.5, p. 263]{Siegmund}
Let  $U\subseteq X, V\subseteq Y$ be open and $f:U\times V \to X$ be such that $f(\cdot,y)$ is integrable for any $y\in V.$ Suppose, that for any $x\in U, y\in V,$ there exist derivative $D_2f(x,y),$ and $$\norm{D_2f(x,y)}\leq g(x),$$
     where $g:U\to \mathbb{R}^+$ is integrable.
    Then the mapping $$F(y) = \int_U f(x,y)dx$$ is differentiable with derivative is given by
    $$DF(y)= \int_UD_2 f(x,y)dx.$$

\end{theorem}

Now let us turn to a key result in random dynamical systems: the fundamental Oseledets's multiplicative ergodic theorem (MET).

In the following we will write $\log^+(x) = \max\{0, \log(x)\}$ for any $x\in X .$ For a probability space $(\Omega, \mathcal{F}, \mathbb{P})$, we denote by $L^1(\Omega, \mathcal{F}, \mathbb{P})$ the space of all integrable measurable functions.

\begin{theorem}[Multiplicative Ergodic Theorem, half line, continuous time]\cite[Thm. 3.4.1, p. 134]{Arnold}\label{MET}

Let $\Phi$ be a linear cocycle  over an ergodic MDS $(\Omega, \mathcal{F}, \mathbb{P}, (\theta_t)_{t \in I})$.	
 Assume that \begin{align}
		\sup_{0\leq t\leq 1}\log^+\norm{\Phi(t,\omega)}&\in L^1(\Omega,\cF,\mathbb{P}),&
		\sup_{0\leq t\leq 1}\log^+\norm{\Phi^{-1}(t,\omega)}&\in L^1(\Omega,\cF,\mathbb{P}).
	\end{align}
 Then there exist a forward invariant set $\Tilde{\Omega}\in\cF$ of full measure such that for every $\omega\in\Tilde{\Omega}$ the following statements hold:
 \begin{enumerate}
    \item[(i)] There exist $k(\omega)$ numbers $\lambda_1(\omega) >...> \lambda_{k(\omega)}(\omega)$ and the invariant splitting $X=U_1(\omega)\oplus...\oplus U_{k(\omega)}(\omega)$ such that
    \begin{align*}
        k(\theta_t\omega) &{=} k(\omega),\\
       \lambda_i(\theta_t\omega) &{=} \lambda_i(\omega) \fall i\in\{1,2,...,k(\omega)\},\\
       d_i(\theta_t\omega) &{=} d_i(\omega) \fall i\in\{1,2,...,k(\omega)\},
    \end{align*}
where $d_i(\omega):=\dim U_i(\omega).$
    \item[(ii)] Put $V_{k(\omega)+1}(\omega):=0$ and for $i=1,..., k(\omega)$
    $$V_i(\omega):=U_i(\omega)\oplus...\oplus U_{k(\omega)}(\omega)$$ so that
    $$V_{k(\omega)}(\omega)\subset ...\subset V_{i}(\omega)\subset ...\subset V_{1}(\omega)$$
    defines a filtration on $X.$ Then for each $x\in X\setminus \{0\}$ the Lyapunov exponent
    $$\lambda(\omega,x) := \lim_{t\to\infty}\frac{1}{t}\log\norm{\Phi(t,\omega)x}$$
    exists and
    $$\lambda(\omega,x) = \lambda_i(\omega)\iff x\in V_i(\omega)\setminus V_{i+1}(\omega),$$
    i.e.$$V_i(\omega) = \{x\in X:\lambda(\omega,x) \leq \lambda_i(\omega)\}.$$
    \item[(iii)] For all $x\in X\setminus \{0\}$
    $$\lambda(\theta_t\omega,\Phi(t,\omega)x) = \lambda(\omega,x),$$
    whence
    $$\Phi(t,\omega) V_i(\omega)\subset V_i(\theta_t\omega)\fall i\in\{1,2,...,k(\omega)\}.$$
 \end{enumerate}
 \end{theorem}

\section*{Acknowledgements}

This work is supported by the Austrian Science Fund
(FWF): DOC 78.

\bibliographystyle{abbrv}
\bibliography{arxiv}

\end{document}